\documentclass[10pt]{amsart}
\usepackage{latexsym}
\usepackage{amsfonts}

\usepackage{amsmath,amsthm,amssymb}

\newtheorem{theor}{Theorem}

\newtheorem{thm}{Theorem}[section]
\newtheorem{lem}[thm]{Lemma}
\newtheorem{cor}[thm]{Corollary}

\newtheorem{prop}[thm]{Proposition}
\theoremstyle{definition}
\newtheorem{defn}[thm]{Definition}

\newtheorem{conv}[thm]{Convention}
\newtheorem{rem}[thm]{Remark}

\newtheorem{propdfn}[thm]{Proposition-Definition}

\begin{document}

\author[Ilya Kapovich]{Ilya Kapovich}

\address{\tt Department of Mathematics, University of Illinois at
  Urbana-Champaign, 1409 West Green Street, Urbana, IL 61801, USA
  \newline http://www.math.uiuc.edu/\~{}kapovich/} \email{\tt
  kapovich@math.uiuc.edu}

\author[Tatiana Nagnibeda]{Tatiana Nagnibeda}

\address{\tt
Section de math\'ematiques,
Universit\'e de Gen\`eve,
2-4, rue du Li\`evre, c.p. 64,
1211 Gen\`eve, Switzerland
\newline http://www.unige.ch/math/folks/nagnibeda}
\email{\tt tatiana.smirnova-nagnibeda@math.unige.ch}

\thanks{The first author was supported by the NSF grant DMS\#0404991.
  Both authors acknowledge the support of the Centre de Recerca
  Matem\`atica at Barcelona and of the Swiss National
  Foundation for Scientific Research}

\title[Patterson-Sullivan embedding]{The Patterson-Sullivan embedding
  and minimal volume entropy for outer space}

\begin{abstract}
  Motivated by Bonahon's result for hyperbolic surfaces, we construct
  an analogue of the Patterson-Sullivan-Bowen-Margulis map from the
  Culler-Vogtmann outer space $CV(F_k)$ into the space of
  projectivized geodesic currents on a free group. We prove that this
  map is a topological embedding and thus obtain a new compactification of the outer space.  
We also prove that for every $k\ge
  2$ the minimum of the volume entropy of the universal covers of
  finite connected volume-one metric graphs with fundamental group of
  rank $k$ and without degree-one vertices is equal to $(3k-3)\log 2$
  and that this minimum is realized by trivalent graphs with all edges
  of equal lengths, and only by such graphs.
\end{abstract}

\subjclass[2000]{ Primary 20F65, Secondary 05C, 37A, 37E, 57M}

\keywords{free groups, metric graphs, Patterson-Sullivan measures,
geodesic currents, volume entropy}

\maketitle

%\tableofcontents

\section{Introduction}\label{intro}

A \emph{geodesic current} on a word-hyperbolic group $G$ is a positive
$G$-invariant Radon measure on the space $\partial^2 G:=\{(x,y):
x,y\in \partial G, x\ne y\}$, where $\partial G$ is the hyperbolic
boundary of $G$ endowed with the canonical boundary topology. The
study of geodesic currents on free groups is motivated by
investigating geometry and dynamics of individual automorphisms, as
well as of groups of automorphisms of a free group. A similar program
proved to be successful in the case of fundamental groups of
hyperbolic surfaces.  Bonahon's foundational work~\cite{Bo86,Bo88}
showed the relevance of geodesic currents to the study of the geometry
of the Teichm\"uller space and of the dynamical properties of surface
homeomorphisms.  Results about geodesic currents in the hyperbolic
surface case can be also found in \cite{Bridg,BT,Ha,Sa1,Pi} and
other sources. Interesting applications of geodesic currents to the
study of free group automorphisms were recently obtained
in~\cite{Ka,Ka1,KLSS,KKS}.

Patterson-Sullivan measures were introduced by Patterson~\cite{Pa76}
and Sullivan~\cite{Su79} in the context of a Kleinian group acting on
the boundary of a hyperbolic space.  The notion was extended by
Coornaert \cite{Coor} to the case of a group $G$ acting geometrically
(that is isometrically, properly discontinuously and cocompactly) on a
Gromov-hyperbolic geodesic metric space (see also related works
\cite{Bou93,BM96,HP,Pau,CP99}). Patterson-Sullivan measures on a
Gromov-hyperbolic space were further studies by Furman~\cite{Fur}.
Patterson-Sullivan measures on the universal covers of finite
simplicial graphs were considered by Lyons~\cite{Ly} and by Coornaert
and Papadopoulos~\cite{CP97}.

Let us remind here briefly the definition in the case of a
non-elementary group $G$ acting geometrically on a $CAT(-1)$ space $X$
(see also definition \ref{defn:cd}.)  For $s>0$ an \emph{$s$-conformal
  density} is a $G$-equivariant family of regular Borel measures
$(\mu_x)_{x\in X}$ on $\partial X$ that are pairwise absolutely
continuous and with the property that their mutual Radon-Nikodym
derivatives satisfy
\[
\frac{d\mu_x}{d \mu_y}(\xi)=e^{-s B_\xi(x,y)}, \quad \text{ for every
} x,y\in X,
\]
where for a point $\xi\in \partial X$ and for $x,y\in X$, $B_\xi(x,y)$
is a Busemann function defined by
\[
B_{\xi}(x,y):=\lim_{z\to\xi, z\in X} [d(x,z)-d(y,z)].
\]
It turns out that there is a unique $s>0$ called the \emph{critical
  exponent} $h(X)$ (see definition \ref{crex}) such that a nonzero
$s$-conformal density exists (and is moreover unique up to scalar
multiplication.) Such a family $(\mu_x)_{x\in X}$ is said to be a
family of \emph{Patterson-Sullivan measures} on $\partial X$. The
critical exponent $h(X)$ coincides with the Hausdorff dimension of
$\partial X$.  

Furman~\cite{Fur} proved, in the more general situation
of a geometric action on a Gromov-hyperbolic space, that there is a
unique up to scalar multiple nonzero $G$-invariant measure $\nu$ on
$\partial^2 X:= \{(\zeta_1,\zeta_2)| \zeta_1,\zeta_2\in \partial X,
\zeta_1\ne \zeta_2\}$ in the same measure class as $\mu_x^2$. Any of
the nonzero scalar multiples of $\nu$ is called an
\emph{$X$-Patterson-Sullivan current.}  Via the identification between
$\partial^2G$ and $\partial^2X$, this measure $\nu$ pulls back to a
canonical, up to a scalar multiple, geodesic current on $G$, any
nonzero scalar multiple of which is called a
\emph{$G$-Patterson-Sullivan current.}

In the case of closed hyperbolic surfaces Patterson-Sullivan currents
admit several other equivalent characterizations.  Let $S$ be a closed
surface with a fixed hyperbolic metric $\rho$, so that
$\widetilde{(S,\rho)}=\mathbb H^2$. Thus $G=\pi_1(S)$ acts on $\mathbb
H^2$ geometrically and $\mathbb H^2/G=S$.  In this situation there is
a natural identification between the space of $G$-invariant measures
on $\partial^2 \mathbb H^2$ and the space of shift-invariant measures
on the unit tangent bundle $\mathbb US$, where the $\mathbb R$-shift
action is given by the geodesic flow on $(S,\rho)$. As shown by
Kaimanovich~\cite{Kaim91}, under this identification
Patterson-Sullivan currents correspond precisely to
\emph{Bowen-Margulis measures} (or \emph{maximal entropy measures}) on
$\mathbb US$, that is the only shift-invariant measures on $\mathbb
US$ whose entropy is equal to the topological entropy of the geodesic
flow on $(S,\rho)$.

For closed hyperbolic surfaces $\mathbb H^2$-Patterson-Sullivan
currents coincide with \emph{Liouville currents} corresponding to the
hyperbolic structure $\rho$.  Bonahon~\cite{Bo86,Bo88} proved that the
map sending a marked hyperbolic structure to the corresponding
projective class of Liouville currents provides a topological
embedding $L: \mathcal T(S)\to \mathbb P Curr(G)$ of the Teichm\"uller
space $\mathcal T(S)$ to the compact space $\mathbb P Curr(G)$ of
projectivized geodesic $G$-currents.

The Culler-Vogtmann outer space~\cite{CV} is a free group analogue of
the Teichm\"uller space. For a free group $F$ of finite rank $k\ge 2$
the \emph{outer space $CV(F)$} consists of equivalence classes of
free, discrete and minimal isometric actions of $F$ on $\mathbb
R$-trees for which the quotient metric graph has volume one.  Two such
actions are equivalent if there is an $F$-equivariant isometry between
the two trees in question.

Let $\Gamma$ be a finite connected graph with no degree-one and
degree-two vertices, and let $\alpha: F\to \pi_1(\Gamma,p)$ be an
isomorphism. Thus $\alpha$ defines a free and discrete action of $F$
on $\tilde\Gamma$ with the quotient $\Gamma$.  Every choice $\mathcal
L$ of a volume-one metric graph structure on $\Gamma$ (that is,
assignment of positive lengths to non-oriented edges of $\Gamma$, so
that the sum of the lengths of all edges is equal to $1$) turns
$\widetilde \Gamma$ into a metric tree so that the above action of $F$
on $\widetilde \Gamma$ becomes an action by isometries.  Hence
$\mathcal L$ defines a point in $CV(F)$. Varying the lengths of edges
of $\Gamma$ gives an open simplex $W_\alpha$ in $CV(F)$ of dimension
$N-1$, where $N$ is the number of non-oriented edges of $\Gamma$.
Thus the outer space $CV(F)$ is a union of open simplices of bounded
dimension.

There is a natural map $\tau: CV(F)\to \mathbb P Curr(F)$ that takes a
point of $CV(F)$ represented by the action of $F$ on a tree, to the
projective class of $F$-Patterson-Sullivan currents corresponding to
this action. We call $\tau$ the \emph{Patterson-Sullivan map} (see Definition \ref{def:PSHmap} for details and for the definition of the Hausdorff dimension map.)  Our
main result is the following statement, which parallels the above
mentioned theorem of Bonahon for hyperbolic surfaces:

\begin{theor}\label{A}
  The Patterson-Sullivan map $\tau: CV(F)\to \mathbb P Curr(F)$ is a
  topological embedding. The Hausdorff dimension map
  $h:CV(F)\to\mathbb R$ is continuous and, moreover, the restriction
  of $h$ to any open simplex in $CV(F)$ is real-analytic.
\end{theor}

Injectivity of $\tau$ follows from a general result of
Furman~\cite{Fur} proved in the context of geometric actions on
Gromov-hyperbolic spaces. The main work in the present paper is in
proving the continuity of $\tau$.

In the case of a closed hyperbolic surface $S$ with $G=\pi_1(S)$
Bonahon proved that the Liouville map $L: \mathcal T(S)\to \mathbb P
Curr(G)$ extends to a homeomorphism from Thurston's compactification
$\widehat{\mathcal T(S)}$ of $\mathcal T(S)$ to the closure of the
image of $L$. It is well-known that $\widehat{\mathcal T(S)}$
coincides with the length-function compactification of $\mathcal
T(S)$.

The map $\tau: CV(F)\to \mathbb P Curr(F)$ is easily seen to be
$Out(F)$-equivariant and the closure of the image of this map
$\overline{\tau(CV(F))}$ is compact. It is therefore natural to ask if
$\tau$ extends to a continuous homeomorphism (that has to be
$Out(F)$-equivariant) from the length function compactification
$\widehat{CV(F)}$ of $CV(F)$ to $\overline{\tau(CV(F))}$. It turns out
that the answer is negative in a very strong sense. Thus Kapovich and
Lustig (in preparation) recently proved that there does not exist a
continuous $Out(F)$-equivariant map $\partial CV(F)\to \mathbb P
Curr(F)$ where $\partial CV(F)=\widehat{CV(F)}-CV(F)$ is the
length-function boundary of $CV(F)$. Therefore any equivariant
topological embedding $CV(F)\to \mathbb P Curr(F)$, such as the
Patterson-Sullivan map $\tau$, results in a new compactification of
$CV(F)$ that is different from the standard length-function
compactification.  This fact is the primary motivation for proving
Theorem~\ref{A} and for obtaining as explicit a description of $\tau$
as possible in the process.

%Two main open questions about the new compactification of the outer
%space which we obtain are that of its contractability and its
%dimension. Recall that $\widehat{CV(F)}$ is contractible and
%finite-dimensional (see the survey article \cite{Vog} for references.)

A different family of continuous
$Out(F)$-equivariant embeddings from $CV(F)$ to $\mathbb P Curr(F)$ was constructed by Reiner Martin~\cite{Ma}.
Unlike the Patterson-Sullivan embedding $\tau$, Martin's embeddings
are not based on a natural geometric construction and use an ad-hoc
procedure, where a point of $CV(F)$ is sent to an explicitly defined
infinite linear combination of ``counting'' currents determined by
conjugacy classes of elements of $F$. His construction leads to
infinite dimensional compactifications of $CV(F)$.

\medskip It is well-understood that in a fairly general negatively
curved setting the Hausdorff dimension of the boundary coincides with
the volume entropy. If $(M,g)$ is a closed connected Riemannian
manifold, then the \emph{volume entropy of $g$} is defined as
\[
h(g):=\liminf_{R\to\infty} \frac{\log {Vol}_{\tilde g}(B(x,R))}{R},
\]
where $x\in \widetilde M$ is a base-point and $B(x,R)$ is the ball of
radius $R$ and center $x$ in $\widetilde M$, equipped with the
pullback $\tilde g$ of the Riemannian metric $g$. It is well known
that the $\liminf$ in the formula can be replaced by $\lim$.  This
definition does not depend on the choice of $x\in \widetilde M$, and
$h(g)>0$ if and only if the group $\pi_1(M)$ has exponential growth.
If $g$ has strictly negative sectional curvature, then $(\widetilde M,
\tilde g)$ is a $CAT(-1)$ space and the Hausdorff dimension of its
boundary (which is also equal to the critical exponent of $(\widetilde
M, \tilde g)$) is equal to the volume entropy $h(g)$. A similar
statement is true for the universal cover of a compact locally
$CAT(-1)$-space $K$. In that case volume has to be interpreted as
counting the number of $\pi_1(K)$-orbit points in the ball of radius
$n$ around the basepoint in $\widetilde K$.

For a compact connected Riemannian manifold $M$ it is natural to ask
what the infimum of $h(g)$ is when $g$ varies over metrics with
$Vol_g(M)=1$ and whether this infimum is achieved.  This is known as
the {\it minimal entropy problem} (see discussion in \cite{BCG}).  A
famous theorem of Besson, Courtois and Gallot~\cite{BCG} shows that if
$M$ admits a locally symmetric volume-one metric $g_0$ of negative
curvature, then $g_0$ minimizes volume entropy among all volume-one
metrics (see an earlier paper of Katok \cite{Kat} for the case of
surfaces.)

A particular case of their theorem (see \cite{BCG1}, Section 5) says
that if $(M, g_0)$ and $(M', g)$ are homotopically equivalent
negatively curved compact connected Riemannian manifolds of the same
dimension $n\geq 3$, and if $(M,g_0)$ is locally symmetric, then
$h^n(g)Vol(M',g)\geq h^n(g_0)Vol(M,g_0)$.  Besson, Courtois and Gallot
also show that $h(g)= h(g_0)$ and $Vol(M',g)=Vol(M,g_0)$ if and only
if $(M',g)$ is isometric to $(M,g_0)$.

In the last section of our paper we prove an analogue of these
statements in the outer space setting.

Theorem~\ref{A} implies that the volume entropy function $h$ (which
again coincides with the Hausdorff dimension of the boundary and with
the critical exponent) factors to a continuous function on the moduli
space $\mathcal M=CV(F)/Out(F)$
\[
\overline h: \mathcal M\to \mathbb R_{>0}.
\]
A point in $\mathcal M$ is a finite connected graph $\Gamma$ without
degree-one and degree-two vertices and with $\pi_1(\Gamma)\cong F$,
endowed with the structure $\mathcal L$ of a volume-one metric graph.
Then $\overline h(\mathcal L)$ is the volume entropy of the metric
tree $\partial \widetilde \Gamma$, where the metric on $\widetilde
\Gamma$ is given by the lift of $\mathcal L$. The analogue of a
locally symmetric manifold is a regular graph (i.e., such that all
vertices are of the same degree) with all edges of equal length.  The
volume entropy of a regular tree with all edges of the same length is
easy to compute explicitly. In particular, assigning the length
$1/(3k-3)$ to each of $(3k-3)$ non-oriented edges in a trivalent graph
with the fundamental group free of rank $k$ gives a volume-one metric
graph with volume entropy of its universal cover equal to $(3k-3)\log
2$.  We prove that this is precisely the minimum of the volume entropy
over all finite connected metric volume-one graphs without vertices of
degree one or two and with fundamental group free of rank $k$.

\begin{rem} Note that while the entropy function is constant on $\mathcal T(S)$, due to the constant curvature, it is not constant on 
  $CV(F)$. One could however make the entropy constant by changing the
  standard normalization (graphs being of total volume $1$) to the
  (less natural) one in which the total volume of a finite connected
  metric graph $\Gamma$ is equal to $1/h(\tilde\Gamma)$.
\end{rem}

\begin{theor}\label{C}
  Let $F$ be a free group of rank $k\ge 2$. Then:
\begin{enumerate}
\item For the function $\overline h: CV(F)/Out(F)\to\mathbb R$ we have
\[
\min \overline h= (3k-3)\log 2.
\]
This minimum is realized by any regular trivalent connected graph
$\Gamma$ with $\pi_1(\Gamma)\cong F$, (so that $\Gamma$ has $3k-3$
non-oriented edges), where each edge of $\Gamma$ is given length
$1/(3k-3)$.

\item If a point of $CV(F)/Out(F)$ realizes the minimum of $h$ then
  this point is a regular trivalent graph with all edges of equal
  lengths.
  
\item We have
\[
\sup_{\mathcal M} \overline h = \infty.
\]
\end{enumerate}
\end{theor}

As an intermediate step in proving Theorem~\ref{C} we establish that
among all the volume-one metric structures on an $m$-regular graph
$\Gamma$ with $m\ge 3$, the volume entropy is minimized by assigning
all the edges of $\Gamma$ equal lengths. This fact is a particular
case of solving the minimal volume entropy problem for an arbitrary
finite connected graph which was carried in a recent work of
Lim~\cite{Lim}.  Related results have also been obtained by
Rivin~\cite{Riv}. The results of Rivin~\cite{Riv} and of
Lim~\cite{Lim} can be used to provide alternative proofs of parts (1)
and (2) of Theorem~\ref{C}. Nonetheless, we choose to present our
proof for completeness and because it uses rather different arguments
from those of Rivin and of Lim (see Section 9.) For the benefit of our topological "outer space" audience we have also incuded an account of Patterson-Sullivan measures (see Section 3), as well as tried to give concrete and explicit proofs to our main
results (Sections 7,8,9.)

The authors are grateful to Florent Balacheff, David Berg, Pierre de
la Harpe, Vadim Kaimanovich, J\'er\^ome Los, Martin Lustig, Paul
Schupp and Dylan Thurston for helpful discussions.  We are also
grateful to Frederic Paulin and Seonhee Lim for informing us about
Seonhee Lim's work and to Jean-Franc\c{o}is Lafont for informing us
about the entropy results of Rivin.

\section{Geodesic currents}

\begin{conv}
  For the remainder of the paper let $F$ be a finitely generated free
  group of rank $k\ge 2$. We will denote by $\partial F$ the space of
  ends of $F$ with the standard ends-space topology. Thus $\partial F$
  is a topological space homeomorphic to the Cantor set. We shall also
  think about $\partial F$ as the hyperbolic boundary of $F$, endowed
  with the canonical boundary topology, in the sense of the theory of
  word-hyperbolic groups (see, for example~\cite{GH}).

  We set
\[
\partial^2 F:=\{(\zeta,\xi): \zeta,\xi\in \partial F \text{ and }
\zeta\ne \xi\}.
\]
\end{conv}

\begin{defn}[Geodesic currents]
  A \emph{geodesic current} on $F$ is a positive, finite on compact
  subsets $F$-invariant Borel measure on $\partial^2 F$.  We denote
  the space of all geodesic currents on $F$ by $Curr(F)$.  The space
  $Curr(F)$ comes equipped with the weak-$*$ topology: for $\nu_n,
  \nu\in Curr(F)$ we have $\displaystyle\lim_{n\to\infty}\nu_n=\nu$
  iff for every two disjoint open sets $S,S'\subseteq \partial F$ we
  have $\displaystyle\lim_{n\to\infty}\nu_n(S\times S')=\nu(S\times
  S')$.
\end{defn}

We say that two nonzero geodesic currents are equivalent, denoted
$\nu_1\sim \nu_2$, if there exists a positive scalar $r\in \mathbb R$
such that $\nu_2=r \nu_1$. We consider also the space \[\mathbb
PCurr(F):=\{\nu\in Curr(F): \nu\ne 0\}/\sim\] \emph{of projectivized
  geodesic currents on $F$}, endowed with the quotient topology.  We
denote the $\sim$-equivalence class of a nonzero geodesic current
$\nu$ by $[\nu]$.

For a (finite or infinite) graph $\Delta$, denote by $V\Delta$ the set
of all vertices of $\Delta$, and denote by $E\Delta$ the set of all
oriented edges of $\Delta$ (i.e., the set of all ordered pairs $(u,v)$
where $u$ and $v$ are adjacent vertices in $\Delta$.)  A path $\gamma$
in $\Delta$ is a sequence of oriented edges which connects a vertex
$o(\gamma)$ (origin) with a vertex $t(\gamma)$ (terminus). A path is
called \emph{reduced} if it does not contain a back-tracking (a path
of the form $(ee^{-1})$).
%The length $|\gamma|$ of the path $\gamma$ is equal
%to the sum of the lengths of all its edges.
We denote by $\mathcal P(\Delta)$ the set of all finite reduced paths
in $\Delta$.
%For an integer $m\ge 1$ we denote by
% $\mathcal P_m(\Delta)$ the set of all reduced paths of
%  length $m$ in $\Delta$.
For a vertex $x\in V\Delta$, we denote by $\mathcal P_x(\Delta)$ the
collection of all $\gamma\in \mathcal P(\Delta)$ that begin with $x$.
For $\gamma\in \mathcal P(\Delta)$, we denote by $a(\gamma)$ the set
of all $e\in E\Delta$ such that $e\gamma\in \mathcal P(\Delta)$ and we
denote by $b(\gamma)$ the set of all $e\in E\Delta$ such that $\gamma
e\in \mathcal P(\Delta)$.

\begin{defn}[Simplicial charts]
  Let $\Gamma$ be a finite connected graph without degree-one vertices
  such that $\pi_1(\Gamma)\cong F$. Let $\alpha:F\to \pi_1(\Gamma,p)$
  be an isomorphism, where $p$ is a vertex of $\Gamma$. We call such
  $\alpha$ a \emph{simplicial chart} for $F$.
\end{defn}

Let $\alpha:F\to \pi_1(\Gamma,p)$ be a simplicial chart.  We consider
$X:=\widetilde \Gamma$, a topological tree, and denote the covering
map from $X$ to $\Gamma$ by $j:X\to \Gamma$. For $\gamma\in \mathcal
P(X)$ we call the reduced path $j(\gamma)$ in $\Gamma$ the
\emph{label} of $\gamma$.  As there is only one reduced path
connecting two arbitrary vertices in a tree, we will often write
$[x,y]$ for a path in $X$ with origin $x$ and terminus $y$.

Let $\partial X$ denote the space of ends of $X$ with the natural
ends-space topology. Then we obtain a canonical $\alpha$-equivariant
homeomorphism $\hat \alpha:\partial F\to\partial X$, as follows.
Suppose we endow $\Gamma$ with the structure of a metric graph, that
is, we assign a positive length to each edge of $\Gamma$. This turns
$X$ into a metric tree and the action of $\pi_1(\Gamma,p)$ on $X$
becomes a discrete isometric action.  Moreover, $X$ is quasi-isometric
to $F$ and, if $F$ is equipped with a word metric and $x_0$ is a lift
of $p$ to $X$, then the orbit map $\tilde\alpha:F\to X$, $f\to
\alpha(f)x_0$, is a quasi-isometry. This quasi-isometry extends to a
homeomorphism $\hat\alpha:\partial F\to \partial X$.  A crucial
feature of this construction is that $\hat\alpha$ does not depend on
the choice of a metric structure on $\Gamma$. If $\alpha$ is fixed, we
will usually suppress explicit mention of $\hat \alpha$ and also of
the map $\alpha$ itself when talking about the action of $F$ on $X$
and on $\partial X$ arising from this situation.  We also denote by
$\partial^2 X$ the set of all pairs $(\zeta_1,\zeta_2)$ such that
$\zeta_1, \zeta_2\in \partial X$ and $\zeta_1\ne \zeta_2$. For
$(\zeta_1,\zeta_2)\in \partial^2 X$ we denote by $[\zeta_1,\zeta_2]$
the simplicial (non-parameterized) geodesic from $\zeta_1$ to
$\zeta_2$ in $X$.  Thus $[\zeta_1,\zeta_2]$ is a subgraph of $X$
isomorphic to the simplicial line, together with a choice of direction
on that line.  We also have the identification $\hat \alpha:
\partial^2 F\to \partial^2 X$.

\begin{defn}[Cylinder sets]
  For every reduced path $\gamma$ in $X$ denote
\begin{gather*}
  Cyl_X(\gamma):= \{(\zeta_1,\zeta_2)\in \partial^2 X: \gamma\subseteq
  [\zeta_1,\zeta_2] \\ \text{ and the orientations on $\gamma$ and on
    $[\zeta_1,\zeta_2]$ agree}\}
\end{gather*}
Also, for $x=o(\gamma)\in X$ denote
\[
Cyl_x(\gamma):=\{\zeta\in \partial X: \gamma \text{ is an initial
  segment of } [x,\zeta]\}
\]
\end{defn}

The collection of all sets $Cyl_X(\gamma)$, where $\gamma$ varies over
$\mathcal P(X)$, gives a basis of closed-open sets for $\partial^2 X$.
For any $x\in X$, the collection of all sets $Cyl_x(\gamma)$, where
$\gamma$ varies over $\mathcal P_x(X)$, gives a basis of closed-open
sets for $\partial X$.  Let us denote
$Cyl_\alpha(\gamma):=\hat\alpha^{-1} Cyl_X(\gamma)$, so that
$Cyl_\alpha(\gamma)\subseteq \partial^2 F$.  It is easy to see that
%\begin{lem}\label{lem:conv}
  %Let $\alpha:F\to \pi_1(\Gamma,p)$ be a simplicial chart and let
  %$X=\widetilde \Gamma$.
for $\nu_n, \nu\in Curr(F)$ $\displaystyle\lim_{n\to\infty}\nu_n=\nu$
iff
$\displaystyle\lim_{n\to\infty}\nu_n(Cyl_\alpha(\gamma))=\nu(Cyl_\alpha(\gamma))$
for every $\gamma \in \mathcal P(X)$.  Moreover, for $\nu,\nu'\in
Curr(F)$ we have $\nu=\nu'$ iff
$\nu(Cyl_\alpha(\gamma))=\nu'(Cyl_\alpha(\gamma))$ for every
$\gamma\in \mathcal P(X)$.
%\end{lem}

\begin{rem}
  Note that for any $f\in F$ and $\gamma\in \mathcal P(X)$ we have $f
  Cyl_\alpha(\gamma)=Cyl_\alpha(f \gamma)$.  Since geodesic currents
  are, by definition, $F$-invariant, for a geodesic current $\nu$ and
  for $\gamma\in \mathcal P(X)$ the value $\nu(Cyl_\alpha(\gamma))$
  only depends on the label $j(\gamma)$ of $\gamma$.
\end{rem}

\section{Patterson-Sullivan measures and metric graphs}\label{sect:ps}

\begin{defn}[Metric and semi-metric graph structures] A \emph{quasi-metric structure $\mathcal L$ on
    a (finite or infinite) graph $\Gamma$} is an assignment of a
  \emph{length} $L(e)\ge 0$ to each edge $e\in E\Gamma$ of $\Gamma$.
  The \emph{volume} of $\mathcal L$ is $vol(\mathcal
  L):=\frac{1}{2}\sum_{e\in E\Gamma} L(e)$. A quasi-metric structure
  is called a \emph{semi-metric structure} if $L(e)=L(e^{-1})$ for
  every $e\in E\Gamma$. A semi-metric structure is called a
  \emph{metric structure} if $L(e)>0$ for every $e\in E\Gamma$.  We
  say that a quasi-metric structure $\mathcal L$ is
  \emph{non-singular} if there is a maximal tree $T$ in $\Gamma$ such
  that $L(e)>0$ for every $e\in E(\Gamma-T)$. A quasi-metric structure
  is \emph{positive} if $L(e)>0$ for every $e\in E\Gamma$.

  If $\mathcal L$ is a semi-metric structure on $\Gamma$, let
  $\Gamma'$ be the graph obtained from $\Gamma$ by contracting to
  points all edges of $\Gamma$ of $\mathcal L$-length zero. Then
  $\Gamma'$ comes equipped with a canonical metric graph structure
  $\mathcal {L'}$ coming from $\mathcal L$. We call $(\Gamma',\mathcal
  {L'})$ \emph{the metric graph associated to $(\Gamma,\mathcal L)$}.
\end{defn}

\begin{conv}\label{conv:semi} Let $\mathcal L$ be a nonsingular semi-metric graph structure
  on a finite graph $\Gamma$. Let $(\Gamma',\mathcal{L'})$ be the
  metric graph associated to $(\Gamma,\mathcal L)$ and let
  $q:\Gamma\to \Gamma'$ be the canonical projection map.
  
Let $X=\widetilde \Gamma$ and let $j:X\to \Gamma$ be the covering
  map. Then $\mathcal L$ lifts canonically to a semi-metric graph structure
  $\widetilde {\mathcal L}$ on $X$ defined as $\widetilde
  {L}(e):=L(j(e))$ for every $e\in EX$. Similarly let $X'=\widetilde
  \Gamma'$ and let $j':X'\to \Gamma'$ be the associated covering map.
  Again, $\mathcal L'$ lifts to a metric graph structure $\widetilde
  {\mathcal L'}$ on $X'$.
  
  It is easy to see that both $j$ and $j'$ preserve edge-lengths and
  that $X'$ is obtained from $X$ by contracting all edges of length
  zero in $X$ to points. Thus $(X', \widetilde{\mathcal L'})$ is the
  metric graph associated to $(X,\widetilde{\mathcal L})$. We denote
  by $\widetilde q: X\to X'$ the canonical projection map.
  
  The semi-metric structure $\widetilde{\mathcal L}$ defines a
  semi-metric $d=d_{\mathcal L}$ on $X$ and $\widetilde{\mathcal L'}$
  defines a metric $d'=d'_{\mathcal L}$ on $X'$. Moreover, $\widetilde
  q: (X,d)\to (X',d')$ is distance-preserving. Note that for both
  $(X,d)$ and $(X',d')$ there are obvious notions of \emph{geodesic
    edge-paths}. In both cases we can metrize $\partial X$ and
  $\partial X'$ by setting
\[
d_x(\xi,\zeta):=e^{-d(x,[\xi,\zeta])} \quad \text{ where }
\xi,\zeta\in \partial X
\]
\[
d'_{x'}(\xi',\zeta'):=e^{-d'(x',[\xi',\zeta'])} \quad \text{ where }
\xi',\zeta'\in \partial X',
\]
where $x\in X,x'\in X'$.  Note that $d_x$ is a metric on $\partial X$,
although $\mathcal L$ was just a semi-metric structure. Moreover, if
$x'=\widetilde q(x)$ then the map $\widetilde q: (\partial X, d_x)\to
(\partial X',d'_{x'})$ is a homeomorphism and an isometry.
\end{conv}

\begin{conv}\label{conv:space}
  For the remainder of this section we will fix $G,X$ and the
  notations below to be one of the following:

\begin{enumerate}
\item We consider a finitely generated group $G$ acting geometrically
  on a $CAT(-1)$ space $X$.
  
  If $x\in X$ is a base-point, the boundary $\partial X$ is metrized
  as follows: for two points $\xi,\zeta\in \partial X$ put
\[
d_x(\xi,\zeta)=
\begin{cases}
  &0, \quad \text{ if } \xi=\zeta,\\
  &exp(-d(x,[\xi,\zeta])), \quad \text{ if } \xi\ne \zeta.
\end{cases}
\]

\item We consider $G=F$ a free group of finite rank $k\ge 2$ and
  $\alpha:F\to \pi_1(\Gamma,p)$ a simplicial chart for $F$, as well as
  a non-singular semi-metric structure $\mathcal L$ defining a
  semi-metric $d$ on $X=\widetilde \Gamma$. Thus $F$ acts on $X$ via
  $\alpha$ by $d$-preserving transformations. In this case let
  $q:\Gamma\to\Gamma', \widetilde q:X\to X',\mathcal{L'}$, $d'$ and
  the metrics on $\partial X$ and $\partial X'$ be as in
  Convention~\ref{conv:semi}. Thus $q_{\#}\circ \alpha: F\to
  \pi_1(\Gamma', p')$ is another simplicial chart for $F$, where
  $p'=q(p)$ and the map $\widetilde q:X\to X'$ is $F$-equivariant.
\end{enumerate}
\end{conv}

Recall from Introduction that, for a point $\xi\in \partial X$ and for
$x,y\in X$, the {\it Busemann function} is
\[
B_{\xi}(x,y):=\lim_{X\ni z\to\xi} \big (d(x,z)-d(y,z)\big ).
\]
If $x,y\in X$ and $\xi\in \partial X$ are such that $y\in [x,\xi]$ then $B_{\xi}(x,y)=d(x,y)$.

Let us denote by $M(\partial X)$ the space of all positive regular
Borel measures on $\partial X$ endowed with the weak-$*$ topology.  If
$\mu\in M(\partial X)$ and $g\in G$ then the measure $g_\ast \mu\in
M(\partial X)$ on $\partial X$ is defined as $(g_\ast
\mu)(A)=\mu(g^{-1}A)$ for a Borel subset $A\subseteq \partial X$.

\begin{propdfn}[Critical Exponent]\label{crex}
  The \emph{Poincar\'e series} of $X$ with respect to a base-point
  $x\in X$ is
\[
\Pi_x(s):=\sum_{g\in G} e^{-s d(x,gx)}\ .
\]
For every $x\in X$ there exists a unique number $h\ge 0$ such that
$\Pi_x(s)$ converges for all $s>h$ and diverges for all $s<h$.  This
number $h$ does not depend on $x\in X$ and is called the
\emph{critical exponent}. We denote it by $h=h(X)=h(G,X)$.
\end{propdfn}

\begin{rem}
  Coornaert discusses this definition in \cite{Coor}. He shows in
  particular that under assumptions of Convention~\ref{conv:space}
  $\Pi_x(h)$ diverges for every $x\in X$, and that the critical exponent
 coincides with the {\it volume entropy} of $X$ defined by the
 right-hand side of the equality:
\[
h(X)=\lim_{R\to\infty} \frac{1}{R} \log \# \{g\in G: d(x,gx)\le R\}
\]
\end{rem}

\begin{defn}[Conformal density]\label{defn:cd}
  For $s\ge 0$, a continuous map $X\to M(\partial X)$, $x\mapsto
  \mu_x$ is called an \emph{$s$-dimensional conformal density on
    $\partial X$ for $G$} if:
\begin{enumerate}
\item The family $(\mu_x)_x$ is $G$-equivariant, that is
  $\mu_{gx}=(g^{-1})_\ast \mu_x$ for every $x\in X,g\in G$.
\item For every $x,y\in X$
\[
\frac{d\mu_x}{d \mu_y}(\xi)=e^{-s B_\xi(x,y)}\ .
\]

\item We have $\mu_x=\mu_y$ if $d(x,y)=0$.
\end{enumerate}
\end{defn}
In particular, we see that for each $x,y\in X$ the measures $\mu_x,
\mu_y$ are absolutely continuous with respect to each other with
bounded Radon-Nikodym derivatives.

The following two statements follow from the basic results established
in \cite{Coor, CP97,BM96}.

\begin{propdfn}[Patterson-Sullivan measures]\label{propdfn:ps}
  The critical exponent $s=h(X)$ is the only value of $s\ge 0$ such
  that there exists a nonzero $s$-dimensional conformal density on
  $\partial X$. Moreover, up to scalar multiplication, the nonzero
  $h$-dimensional conformal density $(\mu_x)_{x}$ is unique.  The
  measures $(\mu_x)_x$ are called \emph{Patterson-Sullivan measures}
  on $\partial X$.
\end{propdfn}

\begin{prop}\label{prop:coor}
  Let $(\mu_x)_x$ be a family of Patterson-Sullivan measures on
  $\partial X$.  Then

\begin{enumerate}
\item The measures $\mu_x$ belong to the same measure class for all
  $x\in X$. Each $\mu_x$ has no atoms and has full support on $\partial
  X$. 
\item For every $x\in X$ the critical exponent $h$ is equal to the
  Hausdorff dimension of $(\partial X,d_x)$. In particular,
  $0<h(X)<\infty$.
\item Let $x,y\in X$ and let $m_y$ be the $h$-dimensional Hausdorff
  measure on $(\partial X,d_y)$. Then $m_y$ and $\mu_x$ are absolutely
  continuous with respect to each other and their mutual Radon-Nikodym
  derivatives are bounded.

%\item For every nonempty open subset $A\subseteq \partial X$ and for
% every $x\in X$ we have $\mu_x(A)>0$ and $m_x(A)>0$.

\end{enumerate}
\end{prop}

Here is another useful characterization of Patterson-Sullivan measures
(see, for example~\cite{Fur}):

\begin{prop}
  Let  $h=h(X)$ be
  the critical exponent, and let $(\mu_x)_x$ be a family of
  Patterson-Sullivan measures on $\partial X$. Then for every $x\in X$
  the measure $\mu_x$ is, up to a scalar multiple, the weak limit as
  $s\to h+$, of the probability measures
  $$\frac{1}{\Pi_x(s)} \sum_{g\in G} e^{-sd(x,gx)} {\rm Dirac}(gx)\ 
  .$$
\end{prop}

\begin{conv} Let us now concentrate our attention on the case when the
acting group is a nonabelian free group $F$ of finite rank $k\geq 2$.
For the remainder of this section, we assume consequently that $F$,
$\Gamma,$ $\Gamma'$, $X$, $X'$, $\mathcal L$ are as in part (2) of
Convention~\ref{conv:space}.

Then
  $h(F,X)=h(F,X')$ since for every $g\in F$ and for every $x\in X$
  with $x'=\widetilde q(x)\in X'$ we have $d(x,gx)=d'(x',gx')$. We
  shall denote this critical exponent by $h(\mathcal L)$.
\end{conv}
  
  Moreover, suppose $(\mu_x)_x$ is a conformal $s$-density on
  $\partial X$. Then for any $x,y\in X$ with $d(x,y)=0$ we have
  $\mu_x=\mu_y$ and hence $(\mu_x)_x$ canonically factors to a
  conformal $s$-density $(\mu'_{x'})_{x'}$ on $\partial X'$.
  Similarly, if $(\mu'_{x'})_{x'}$ is a conformal $s$-density on
  $\partial X'$, then it canonically pulls back to a conformal
  $s$-density $(\mu_x)_x =((\tilde q^{-1})_*\mu'_{q(x)})_x$ on
  $\partial X$.
  
\begin{rem}\label{simp}
  Let $\gamma\in \mathcal P(X)$ be an edge-path from a vertex $x$ to a
  vertex $y$ in $X$. When restricted to $Cyl_x(\gamma)$, condition (2)
  of Definition~\ref{defn:cd} simplifies to
\[
\mu_x=e^{-s d(x,y)}\mu_y \quad \text{ on } Cyl_x(\gamma)
\]
and, in particular,
\[
\mu_x(Cyl_x(\gamma))=e^{-s d(x,y)}\mu_y (Cyl_x(\gamma)).
\]

This shows that an $s$-conformal density $(\mu_x)_x$ is uniquely
determined by the values $\mu_x(Cyl_x(f))$, where $x$ varies over the
vertices of $X$ and $f$ varies over all edges of $X$ with origin $x$.
Moreover, in view of $F$-equivariance of $(\mu_x)_x$, it suffices to
take $x$ from a bijective lift of the vertex set of $\Gamma$ to $X$.
\end{rem}

\begin{conv}\label{conv:measures}
  Let $e$ be an oriented edge of $\Gamma$ and let $(\mu_x)_x$ be an
  $s$-conformal density for $X$ with $s>0$. Let $f$ be a lift of $e$
  to $X$ and let $x$ be the origin of $f$. We denote
  $$w_e=w_{e,\mathcal L}:=\mu_x(Cyl_x(f))\ .$$
  Because of
  $F$-equivariance of $(\mu_x)_x$ the value $w_e$ does not depend on
  the choice of the lift $f$ of $e$.
\end{conv}

\begin{prop}\label{eq}
  Let $s>0$. Let $(\mu_x)_x$ be an $s$-conformal density for $X$.
  Then, for every $e\in E\Gamma$,
\begin{enumerate}
\item we have $w_e>0$;
\item we have
\[
w_e = \exp(-s L(e)) \sum_{e'\in b(e)} w_{e'}. \tag{*}
\]
\end{enumerate}
Moreover, if $(w_e)_{e\in E\Gamma}$ satisfy conditions (1), (2) above,
then there exists a unique $s$-conformal density $(\mu_x)_x$ such that
for every $e\in E\Gamma$ and for every lift $f$ of $e$ to $X$ with
origin $x$ we have $w_e=\mu_x(Cyl_x(f))$.
\end{prop}

\begin{proof}
  Suppose $(\mu_x)_x$ is an $s$-conformal density for $X$. Let $e$ be
  an edge of $\Gamma$ and let $f$ be a lift of $e$ to $X$ with origin
  $x$. Since $Cyl_x(f)\subseteq \partial X$ is a nonempty open set,
  Proposition~\ref{prop:coor} implies that $w_e=\mu_x(Cyl_x(f))>0$ so
  that condition (1) holds.  Let $y$ be the terminal vertex of $f$.
  For every edge $e'\in b(e)$ there is a unique lift $f'$ of $e'$ to
  $X$ with origin $y$.  Then
\[
Cyl_x(f)=\bigsqcup_{f'} Cyl_x(ff')
\]
and hence
\[
\mu_x(Cyl_x(f))=\sum_{f'} \mu_x(Cyl_x(ff')).
\]
But $Cyl_x(ff')=Cyl_y(f')$ and we have
\[
\mu_x(Cyl_x(ff'))=e^{-s d(x,y)}\mu_y(Cyl_y(f'))=e^{-s
  L(e)}\mu_y(Cyl_y(f')).
\]
Therefore
\[
w_e = e^{-s L(e)} \sum_{e'\in b(e)} w_{e'}\ .
\]
If $(w_e)_{e\in E\Gamma}$ satisfy conditions (1) and (2), then it is
not hard to check that using formulae from Remark~\ref{simp} one can
define an $s$-dimensional conformal density $(\mu_x)_{x\in X}$, as
required.  We leave the details of this verification to the reader.
\end{proof}

We conclude this section with a short note on Hausdorff measures.  It
follows from the definitions that for any $x,y\in X$ the metrics
$d_x$, $d_y$ on $\partial X$ are Lipschitz-equivalent to each other
and hence have the same Hausdorff dimension. Let $s>0$ and let
$\mathcal H^s_x$ be the $s$-dimensional Hausdorff measure on
$(\partial X,d_x)$.

Let $\gamma=[x,y]$ be a geodesic segment in $X$. Then $d_x =
e^{-d(x,y)}d_y \quad \text{ on } Cyl_x(\gamma)\subseteq \partial X.$
Therefore, by definition of Hausdorff measures,
\[
\mathcal H_x^s=e^{-sd(x,y)}\mathcal H_y^s \text{ on } Cyl_x(\gamma).
\]
Thus, for $s$ equal to the Hausdorff dimension of $\partial X$, the
family $(\mathcal H^s_x)_{x}$ is a nonzero $s$-dimensional conformal
density and provides a family of Patterson-Sullivan measures on
$\partial X$.  In particular, if we take a lift $f$ of every edge $e$
to $X$ with origin $x$ and denote $\theta_{e,s}:=\mathcal
H_x^s(Cyl(f))$ for $s\geq 0$, then the numbers $\theta_{e,s}$ satisfy
the system of equations \thetag{*} from Proposition~\ref{eq}:
\[
\theta_{e,s}=\sum_{e'\in b(e)}\exp(-sL(e)) \theta_{e',s} \quad e\in
E\Gamma.
\]

\section{Perron-Frobenius theory for metric trees}\label{sect:pf}

Systems of equations of the type \thetag{*} appearing in
Proposition~\ref{eq} arise in various contexts and can be studied by
the theory of Perron-Frobenius-Ruelle.  The matrix $A_\mathcal L(s)$
of such a system (see Convention \ref{conv:matrix} for the precise definition) 
is a transfer operator, and the statements of
Lemma~\ref{md} and Corollary~\ref{cor:md} below are standard facts
about transfer operators (see for example the article of Guillop\'e
\cite{Gu} where dynamics on metric trees is studied in detail.) In the
probabilistic setting, Perron-Frobenius theory can be applied to study
random walks on trees with finitely many cone types (among them universal covers of finite graphs.) 
In particular it allows the computation of the
rate of escape of a random walk and of the spectral radius of its
transition operator, see \cite{N, NW}.

Below we shall give a self-contained exposition of the basic facts
from the Perron-Frobenius theory that we need (see \cite{Se} for a
detailed exposition).  We shall adapt to our situation the approach of
Edgar \cite{E} to the study of self-similar fractals through so-called
Mauldin-Williams graphs~\cite{MW}. In particular, the proof of
Lemma~\ref{md} below follows closely the proof of Theorem 6.6.6 in
\cite{E}.

If $A$ is a matrix with real coefficients, we will denote by $r(A)$
the \emph{spectral radius} of $A$.  Recall that a nonnegative matrix
$A$ is called \emph{irreducible} if for every position $ij$ there
exists an integer $n>0$ such that $(A^n)_{ij}>0$.  
The notation $A\ge 0$ means that all entries of $A$ are nonnegative
and the notation $A>0$ means that all entries of $A$ are positive. If
$A$ and $B$ are matrices of the same size, we write $A\le B$ if
$B-A\ge 0$ and $A<B$ if $B-A>0$.

\begin{propdfn}[Perron-Frobenius Theorem]\label{pf}
  Let $A\ge 0$ be an irreducible nonnegative $n\times n$-matrix,
  $n\geq 1$.  Then:

\begin{enumerate}
\item The number $r(A)>0$ is an eigenvalue of $A$ of multiplicity $1$.
\item There exists a (unique up to a scalar multiple) column vector
  $Y>0$ such that $AY=r(A)Y$.
\item If $Y=\begin{bmatrix} y_1 \\ \vdots\\ y_n \end{bmatrix} \ge 0$,
  a nonzero column vector, and $\lambda\in \mathbb R$ are such that
  $AY=\lambda Y$, then $\lambda=r(A)$.
\item Suppose that $Y\ge 0$, a nonzero column vector and $\lambda\in
  \mathbb R$ are such that $AY\le \lambda Y$ and such that for some
  coordinate $i$ we have $(AY)_i< \lambda y_i$. Then $r(A)<\lambda$.
\item Suppose that $Y\ge 0$, a nonzero column vector, and $\lambda\in
  \mathbb R$ are such that $AY\ge \lambda Y$. Then $r(A)\ge \lambda$.
\end{enumerate}

The number $r(A)$ is called the \emph{Perron-Frobenius eigenvalue} of
$A$. A column eigenvector $Y>0$ such that $AY=r(A)Y$ is called a
\emph{(right) Perron-Frobenius eigenvector} of $A$.

In this situation the transposed matrix $A^T$ is also irreducible and
$r(A)=r(A^T)$, so that $A$ and $A^T$ have the same Perron-Frobenius
eigenvalue.  If $U$ is a right Perron-Frobenius eigenvector of $A^T$,
the row-vector $U^T$ is called a \emph{left Perron-Frobenius
  eigenvector} of $A$.
\end{propdfn}

\begin{conv}\label{conv:matrix}
  For the remainder of this section, unless specified otherwise, let
  $F,\Gamma,X$ be as in part (2) of Convention~\ref{conv:space}.  Let
  $n=\#E\Gamma$ be the number of oriented edges of $\Gamma$ and let us
  fix an ordering $e_1,\dots, e_n$ on $E\Gamma$. Also, let $\mathcal
  L$ be a positive quasi-metric or a non-singular semi-metric
  structure on $\Gamma$.
  
  Let $H(\Gamma)$ denote the \emph{reduced line graph} of $\Gamma$,
  which is defined as follows. The vertex set of $H(\Gamma)$ is
  $E\Gamma$. The set of oriented edges of $H(\Gamma)$ consists of
  all reduced paths in $\Gamma$ containing exactly two edges.  An edge
  $\gamma=ee'$ of $H(\Gamma)$ has the origin $e$ and the terminus
  $e'$. The inverse edge of $\gamma$ is the path $(e')^{-1}e^{-1}$.
  Let $M$ be the adjacency matrix of $H(\Gamma)$, that is $M=(m_{ij})_{i,j=1}^n$, where
\[
m_{ij}:=\begin{cases} &1, \quad \text{if } e_ie_j\in \mathcal
  P(\Gamma)\\
  &0,\quad \text{otherwise.}
\end{cases}
\]
Denote $A_{\mathcal L}(s):=Diag(e^{-sL(e_1)}, \dots, e^{-sL(e_n)})M$.
The system \thetag{*} from part (2) of Proposition~\ref{eq} rewrites
as the matrix equation:
\[
A_{\mathcal L}(s) \begin{bmatrix} w_{e_1} \\ \vdots \\
  w_{e_n}\end{bmatrix}=\begin{bmatrix} w_{e_1} \\ \vdots \\
  w_{e_n}\end{bmatrix}.
\]
Let $\Phi_{\mathcal L}(s)$ denote the spectral radius of $A_{\mathcal
  L}(s)$.
\end{conv}

\begin{lem}\label{md}
  The following hold:
\begin{enumerate}
\item The matrices $A_{\mathcal L}(s)$ and $A_{\mathcal L}(s)^T$ are
  nonnegative and irreducible for every $s\in \mathbb R$.
  
\item The function $\Phi_{\mathcal L}(s)$ is continuous and strictly
  monotone decreasing on the interval $0\le s<\infty$.
  
\item We have $\Phi_{\mathcal L}(0)>1$.
  
\item If $\mathcal L$ is a positive quasi-metric structure then
  $\lim_{s\to\infty}\Phi_{\mathcal L}(s)=0$.

\end{enumerate}
\end{lem}
\begin{proof}
  Recall that $\Gamma$ is finite, connected, has no degree-one
  vertices and $\pi_1(\Gamma)$ is a free group of rank $k\ge 2$.
  Therefore the graph $H(\Gamma)$ is strongly connected and hence its
  adjacency matrix $M$ is nonnegative irreducible and the same is true
  for its transpose $M^T$. The matrix $A_{\mathcal L}(s)$ is obtained
  from $M$ by multiplying the $i$-th row of $M$ by a positive number
  $e^{-sL(e_i)}$ for each $i=1,\dots, n$. Hence $A_{\mathcal L}(s)$
  and $A_{\mathcal L}(s)^T$ are nonnegative and irreducible.
  
  The continuity of $\Phi_{\mathcal L}(s)$ follows from its
  definition.
  
  Suppose now that $0\le s<s'$.  Let $Y=\begin{bmatrix} y_1\\ \vdots
    \\ y_n
\end{bmatrix}$ be a positive Perron-Frobenius eigenvector of
$A_{\mathcal L}(s)$, so that $A_{\mathcal L}(s)Y=\Phi_{\mathcal
  L}(s)Y$.  Since $L(e_i)\ge 0$, the functions $e^{-sL(e_i)}$ are
monotone non-increasing for each $i$.  Hence component-wise
$a_{ij}(s)\le a_{ij}(s')$ and therefore $A_{\mathcal L}(s')Y\le
A_{\mathcal L}(s)Y=\Phi_{\mathcal L}(s)Y$.  Moreover, there is some
edge $e_i$ with $L(e_i)>0$ and hence $[A_{\mathcal
  L}(s')Y]_i<[A_{\mathcal L}(s)Y]_i=\Phi_{\mathcal L}(s) y_i$.
Therefore $\Phi_{\mathcal L}(s')<\Phi_{\mathcal L}(s)$, as claimed.

Note that $A_{\mathcal L}(0)=M$ and $\Phi_{\mathcal L}(0)$ is the
Perron-Frobenius eigenvalue of $M$. The fundamental group of
$H(\Gamma)$ is free of rank at least two. Hence the universal cover of
$H(\Gamma)$ has exponential growth, that is, the spectral radius of
$M$ is bigger than $1$.

To see that (4) holds, note that if $L(e)>0$ for every $e\in E\Gamma$
then we have $\displaystyle\lim_{s\to\infty} \sum_{ij} a_{ij}(s)=0$.
Also,
\[A(s)\begin{bmatrix} 1\\ 1\\ \vdots\\ 1\end{bmatrix}\le \sum_{ij}
a_{ij}(s)\begin{bmatrix} 1\\ 1\\ \vdots\\ 1\end{bmatrix}.\] Hence
$\Phi(s)\le \sum_{ij} a_{ij}(s)$ and so $\lim_{s\to\infty} \Phi(s)=0$.
\end{proof}

\begin{cor}\label{cor:md}
  Let $\mathcal L$ be a positive quasi-metric structure or a
  non-singular semi-metric structure on $\Gamma$.  Then there exists a
  unique $s>0$ such that $\Phi_{\mathcal L}(s)=1$. If $\mathcal L$ is
  a non-singular semi-metric structure then $s=h(\mathcal L)$.
\end{cor}
\begin{proof}
  Lemma~\ref{md} implies that there is at most one $s>0$ such that
  $\Phi_{\mathcal L}(s)=1$.  If $\mathcal L$ is a non-singular
  semi-metric structure, the existence of Patterson-Sullivan measures
  (Proposition-Definition~\ref{propdfn:ps}) and Proposition~\ref{eq}
  guarantee that when $s=h(\mathcal L)$, the Perron-Frobenius
  eigenvalue of $A_{\mathcal L}(s)$ is equal to $1$, that is, that
  $\Phi_{\mathcal L}(h(\mathcal L))=1$.  If $\mathcal L$ is a positive
  quasi-metric structure then parts (2), (3) and (4) of Lemma~\ref{md}
  guarantee the existence of $s>0$ such that $\Phi_{\mathcal L}(s)=1$.
\end{proof}

>From now on, given a positive quasi-metric structure (or a
non-singular semi-metric structure) on $\Gamma$ we will denote by
$h(\mathcal L)$ the unique value $s>0$ such that $\Phi_{\mathcal
  L}(s)=1$, and will refer to $h(\mathcal L)$ as the \emph{volume
  entropy} of $\mathcal L$.

\begin{rem} Note that if $\mathcal L$ is the metric structure which assigns the same lengths to all edges in $\Gamma$, this description of $h(\mathcal L)$ specializes to the explicit formula known for the volume entropy of uniform simplicial trees, which comes from the consideration of the corresponding subshift of finite type (see e.g. \cite{Ly}.) 
\end{rem}

We will now rewrite the system  $A_{\mathcal L}(s)Y=Y$ in
the form allowing to apply the Implicit Function Theorem. This system
is equivalent to the following $n$ equations:
\[
e^{-s L(e_i)}(m_{i1}y_1+...+m_{in}y_n)-y_i=0, \quad i=1,\dots, n.
\]
To express $y_1,\dots, y_n,s$ as implicit functions of $L(e_1),\dots,
L(e_n)$ we need an extra normalizing equation: $y_1^2+\dots +y_n^2=1$.

\begin{prop}\label{det}  Let $L_1=L(e_1),\dots, L_n=L(e_n)$ be a non-singular semi-metric
  structure or a positive quasi-metric structure $\mathcal L$ on
  $\Gamma$. We set
\[F_i(L_1,\dots,L_n,y_1,\dots,y_n,s):=e^{-s L_i}(m_{i1}y_1+...+m_{in}y_n)-y_i\] for $i=1,\dots, n$, and
\[F_{n+1}(L_1,\dots,L_n,y_1,\dots,y_n,s):=y_1^2+\dots +y_n^2-1.\]
Consider the following system of $n+1$ equations in $2n+1$ variables:
\[
F_i(L_1,\dots,L_n,y_1,\dots,y_n,s)=0,\quad i=1,\dots, n+1\tag{!}
\]
Let $J$ be the Jacobian of this system, that is the
$(n+1)\times(n+1)$-matrix consisting of the partial derivatives of
$F_1,\dots, F_{n+1}$ with respect to $y_1,\dots, y_n,s$:
\[
J_{ij}=\begin{cases} & \frac{\partial F_i}{\partial y_j}\quad 1\le
  i\le n+1, 1\le j\le n\\
  &\frac{\partial F_i}{\partial s}\quad 1\le i\le n+1, j=n+1.
\end{cases}
\]
Suppose $s>0, y_i>0$, for $i=1,\dots, n$, are such that 
$z=(L_1,\dots,L_n,y_1,\dots,y_n,s)$ satisfies the system \thetag{!}.
Then $\det J|_z\ne 0$.
\end{prop}
\begin{proof}
  Let us compute the matrix $J$ at $z$, using the information that $z$
  satisfies \thetag{!}. We will denote $a_{ij}=(A_{\mathcal
    L}(s))_{ij}=e^{-sL_i}m_{ij}$.
  
  For $i\ne j$, $1\le i,j\le n$ we get $\frac{\partial F_i}{\partial
    y_j}=e^{-sL_i}m_{ij}=a_{ij}$. For $i=j$ we get $\frac{\partial
    F_i}{\partial y_i}=e^{-sL_i}m_{ii}-1=a_{ii}-1$. Thus in the upper
  left corner of $J$ we see the $n\times n$ matrix $A_{\mathcal
    L}(s)-I_n$.

  Let us compute $\frac{\partial F_i}{\partial s}$. We have
\[
\frac{\partial F_i}{\partial s}=-L_ie^{-s
  L_i}(m_{i1}y_1+...+m_{in}y_n)=-L_iy_i\quad \text{ for } i=1,\dots,
n,
\]
where the last equality holds since $F_i(z)=0$.

Finally, the last row of $J$ obtained by differentiating
$F_{n+1}=y_1^2+\dots +y_n^2-1$ along $y_1,..,y_n,s$ is $[2y_1 \ 2y_2 \ 
\dots \ 2y_n\ 0]$.

Thus
\[
J=\begin{bmatrix} a_{11}-1 & a_{12} & a_{13} & \dots & a_{1n}&
  -L_1y_1\\
  a_{21} & a_{22}-1& a_{23}& \dots & a_{2n} & -L_2 y_2\\
  \vdots & \vdots & \vdots & \dots& \vdots& \vdots\\
  a_{i1} & a_{i2} & a_{i3} & \dots & a_{in}& -L_iy_i\\
  \vdots & \vdots & \vdots & \dots& \vdots& \vdots\\
  a_{n1} & a_{n2} & a_{n3} & \dots & a_{nn}-1& -L_ny_n\\
  2y_1 & 2y_2 & 2y_3&\dots & 2y_n& 0
\end{bmatrix}
\]
We claim that the rows of the matrix $J$ are linearly independent and
hence $\det J\ne 0$. The column vector $Y=\begin{bmatrix} y_1\\ \vdots
  \\ y_n
\end{bmatrix}$ satisfies
$(A_{\mathcal L}(s)-I_n)Y=0$. This implies that the last row of $J$ is
perpendicular to the first $n$ rows.

Since $Y>0$, it therefore suffices to show that the first $n$ rows of
$J$ are linearly independent.

Note that $\det(A_{\mathcal L}(s)-I_n)=0$. However the matrix
$A_{\mathcal L}(s)-I_n$ has rank $n-1$ since $1$ is the
Perron-Frobenius eigenvalue of $A_{\mathcal L}(s)$ and hence has
multiplicity one.  Thus, up to a scalar, there is only one nontrivial
linear relation between the rows of $A_{\mathcal L}(s)-I_n$. This
relation is given by the left Perron-Frobenius eigenvector
$Z=[z_1,\dots , z_n]$ of $A_{\mathcal L}(s)$. Indeed $ZA_{\mathcal
  L}(s)=Z$ and $Z[A_{\mathcal L}(s)-I_n]=0$.  Note that $z_i>0$ for
all $i=1,\dots, n$.

Suppose that the first $n$ rows of $J$ are linearly dependent and that
we have a nonzero row vector $Z'$ of length $n$ such that $Z'J_n=0$
where $J_n$ is the $n\times (n+1)$ matrix consisting of the first $n$
rows of $J$. Then $Z'$ is also a relation between the first $n$ rows
of $A_{\mathcal L}(s)-I_n$ and hence $Z'$ is a multiple of $Z$. Thus
$ZJ_n=0$.

However, when we multiply $Z$ by the last column of $J_n$ to compute
the $(n+1)$-st entry in $ZJ_n$, we get $-L_1y_1z_1-L_2y_2z_2-\dots
-L_ny_nz_n$.

This number is strictly negative since $L_i\ge 0, y_i>0, z_i>0$ for
all $i=1,\dots, n$ and there is some $i$ such that $L_i>0$.  This
gives us a contradiction with the fact that $ZJ_n=0$.
\end{proof}

For the remainder of this section we will denote an $n$-tuple
$(p_1,\dots, p_n)\in \mathbb R^n$ by $\overline p$.

\begin{cor}\label{cor:use}
  Let $L_1^{(0)}=L^{(0)}(e_1),\dots, L_n^{(0)}=L^{(0)}(e_n)$ be a
  non-singular semi-metric or a positive quasi-metric structure
  $\mathcal L^{(0)}$ on $\Gamma$.  Suppose $s^{(0)}>0, y_i^{(0)}>0$,
  where $i=1,\dots, n$, are such that the point $z^{(0)}=({\overline
    L}^{(0)},{\overline y}^{(0)},s^{(0)})\in \mathbb R^{2n+1}$
  satisfies the system \thetag{!}. Then there exist an open
  neighborhood $U$ of ${\overline L}^{(0)}$ in $\mathbb R^n$ and
  real-analytic functions $s=s(\overline{L})$, $y_i=y_i(\overline{L})$
  on $U$ such that for every $\overline L\in U$ the point
  \[(\overline{L}, y_1(\overline{L}),
  \dots, y_n(\overline{L}), s(\overline{L}))\in \mathbb R^{2n+1}\]
  satisfies \thetag{!} and such that $y_i({\overline
    L}^{(0)})=y_i^{(0)}$, $s({\overline L}^{(0)})=s^{(0)}$.
  
  Moreover, whenever ${\overline L}\in U$ defines a non-singular
  semi-metric structure $\mathcal L$ on $\Gamma$, $s({\overline
    L})$ is equal to the critical exponent of $(\widetilde \Gamma,
  d_{\mathcal L})$ and $(y_1({\overline L}), \dots, y_n({\overline
    L}))$ is a scalar multiple of $(w_{e_1,\mathcal L}, \dots,
  w_{e_n,\mathcal L})$.
\end{cor}
\begin{proof}
  Proposition~\ref{det} implies that the Implicit Function Theorem is
  applicable at the point $z^{(0)}$. Thus there exists an open
  neighborhood $U$ of $z^{(0)}$ and real-analytic functions
  $s=s({\overline L})$, $y_i=y_i({\overline L})$ on $U$ such that for
  every ${\overline L}\in U$ the point \[({\overline L},
  y_1({\overline L}), \dots, y_n({\overline L}), s({\overline L}))\in
  \mathbb R^{2n+1}\] satisfies (!)  and such that $y_i({\overline
    L}^{(0)})=y_i^{(0)}$, $s({\overline L}^{(0)})=s^{(0)}$.

  Moreover, since $y_i^{(0)}>0$, we can choose $U$ so that
  $y_i=y_i({\overline L})>0$ on $U$. Let ${\overline L}\in U$ define a
  non-singular semi-metric structure $\mathcal L$ on $\Gamma$. By
  Proposition~\ref{eq} the critical exponent $h=h(\mathcal L)$ of
  $(\widetilde \Gamma, d_{\mathcal L})$ satisfies the property that
  $\Phi_{\mathcal L}(h)=1$.  Also, by construction, $\Phi_{\mathcal
    L}(s({\overline L}))=1$.  Corollary~\ref{cor:md} now implies that
  $h=s({\overline L})$.  Moreover, Proposition~\ref{eq} and the
  definition of the functions $y_i(L_1,\dots, L_n)$ imply that both
$(w_{e_1,\mathcal L}, \dots, w_{e_n,\mathcal L})$ and
$(y_1({\overline L}), \dots, y_n({\overline L}))$
are Perron-Frobenius eigenvectors of the matrix $A_{\mathcal L}(h)$.
Therefore they are scalar multiples of each other, as required.
\end{proof}

\section{Patterson-Sullivan currents}

The following is essentially a corollary of Proposition~1 of
Furman~\cite{Fur}.

\begin{propdfn}[Patterson-Sullivan current]
  Let $G,X$ be as in Convention~\ref{conv:space}. Let $(\mu_x)_{x\in
    X}$ be a family of Patterson-Sullivan measures on $\partial X$ and
  let $\mu=\mu_y$ for some $y\in X$. Then there exists a unique, up to
  a scalar multiple, $G$-invariant and flip-invariant nonzero locally
  finite measure $\nu$ on $\partial^2 X$ in the measure class of
  $\mu\times \mu$.
  
  Moreover, this measure $\nu$ is of the form
\[
d\nu(\xi,\zeta)=e^{2h(X) f_{\mu}(\xi,\zeta)} d\mu(\xi) d\mu(\zeta),
\]
where $f_{\mu}:\partial^2 X\to \mathbb R_+$ is a symmetric Borel
function which is within bounded distance from the function $d(x,
[\xi,\zeta])$.  Such a measure $\nu$ is called an
\emph{$X$-Patterson-Sullivan current} for the action of $G$ on $X$.
Since $\nu$ is unique up to a scalar multiple, its projective class
$[\nu]$ is called \emph{the projective $X$-Patterson-Sullivan
  current}.
\end{propdfn}

For the remaining part of this section let $F, \Gamma, \mathcal L, X,
d, \alpha$ be as in Part 2) of Convention~\ref{conv:space}.

\begin{defn}\label{def:Fcurr}
  Recall that the choice of simplicial chart $\alpha$ defines a
  homeomorphism $\hat \alpha: \partial^2 F\to\partial^2 X$. Let $\nu$
  be an $X$-Patterson-Sullivan current. Then its pull-back $\hat
  \alpha_\ast(\nu)$ is an $F$-invariant measure on $\partial^2 F$
  which is called an \emph{$F$-Patterson-Sullivan current} for the
  pair $(\alpha, \mathcal L)$.  Its projective class $[\nu]$ is called
  \emph{the projective $F$-Patterson-Sullivan current} for the pair
  $(\alpha, \mathcal L)$.
\end{defn}

We now proceed to give an explicit formula for the
$X$-Patterson-Sullivan current associated with the action of $F$ on
$X$.

\begin{prop}\label{prop:pscurr}
  Let $z\in X$, and let $h(\mathcal L)$ be the critical exponent of
  $X$. Let $(\mu_x)_x$ be a family of Patterson-Sullivan measures on
  $\partial X$ and let $w_e$ be defined as in
  Convention~\ref{conv:measures}.
  
  Then the measure $\nu$ on $\partial^2 X$ given by the formula
\[
d\nu(\xi,\zeta)=e^{2h(\mathcal L) d(z, [\xi,\zeta])} d\mu_z(\xi)
d\mu_z(\zeta)\tag{$\clubsuit$}
\]
is an $X$-Patterson-Sullivan current.

Moreover, for any path $\gamma=[x,y]\in \mathcal P(X)$  we have
\[
\nu(Cyl_X(\gamma))=e^{-h(\mathcal L) L(\gamma)} \big (\sum_{e\in
  b(e')} w_e\big ) \big (\sum_{e\in b(e'')} w_e\big )\tag{$\dag$}
\]
where $(e')^{-1}\in E\Gamma$ is the label of the first edge of
$\gamma$ and $e''\in E\Gamma$ is the label of the last edge of
$\gamma$.
\end{prop}

\begin{proof}
  We will first show that $\thetag{\dag}$ defines a geodesic current
  on $\partial X$. That is, we claim that there exists a unique
  geodesic current $\nu'$ such that for every $\gamma$ as in the
  statement of the proposition
\[
\nu'(Cyl_X(\gamma))=e^{-h(\mathcal L) L(\gamma)} \big (\sum_{e\in
  b(e')} w_e\big ) \big (\sum_{e\in b(e'')} w_e\big ).
\]
%Recall that $Cyl_X(\gamma)\subseteq \partial^2X$ and
%$Cyl_x(\gamma)\subseteq \partial X$.
In view of the definition of $w_e$'s the above formula is equivalent
to
\[
\nu'(Cyl_X(\gamma))=e^{-h(\mathcal L) L(\gamma)}
\mu_x(Cyl_y(\gamma^{-1})) \mu_y(Cyl_x(\gamma)).\tag{$\ddag$}
\]
The uniqueness of $\nu'$ is obvious. Also, by construction $\nu'$ is
$F$-invariant, provided that $\nu'$ is a measure. Thus it remains to
show that the above formula does define a measure on $\partial^2 X$.
To do this we need to check (see, for example, \cite{Ka1}) that for
every $\gamma$ as above
\[
\nu'(Cyl_X(\gamma))=\sum_{f\in b(\gamma)}\nu'(Cyl_X(\gamma f))
\]
and
\[
\nu'(Cyl_X(\gamma))=\sum_{f\in a(\gamma)}\nu'(Cyl_X(f\gamma)).
\]
We will verify the first formula, as the second one is completely
analogous. By $\thetag{\ddag}$ applied to each of the paths $\gamma
f$, where $f\in b(\gamma)$, we have

\begin{gather*}
  \nu'(Cyl_X(\gamma f))=e^{-h(\mathcal L) L(\gamma f)}
  \mu_x(Cyl_{t(f)}(f^{-1}\gamma^{-1})) \mu_{t(f)}(Cyl_x(\gamma
  f))=\\
  e^{-h(\mathcal L) L(\gamma
    f)}\mu_x(Cyl_{y}(\gamma^{-1}))\mu_{t(f)}(Cyl_x(\gamma f))=\\
  e^{-h(\mathcal L) L(\gamma f)}\mu_x(Cyl_{y}(\gamma^{-1}))
  \mu_{t(f)}(Cyl_y(f)).
\end{gather*}

Since
\[
Cyl_x(\gamma)=\bigsqcup_{f\in b(\gamma)} Cyl_x(\gamma f),
\]
it follows that
\[
\mu_y(Cyl_x(\gamma))=\sum_{f\in b(\gamma)} \mu_y(Cyl_x(\gamma f)).
\]
Therefore
\begin{gather*}
  \nu'(Cyl_X(\gamma))=e^{-h(\mathcal L) L(\gamma)}
  \mu_x(Cyl_y(\gamma^{-1}))
  \mu_y(Cyl_x(\gamma))=\\
  e^{-h(\mathcal L) L(\gamma)} \mu_x(Cyl_y(\gamma^{-1}))\big
  (\sum_{f\in b(\gamma)}
  \mu_y(Cyl_x(\gamma f))\big )=\\
  e^{-h(\mathcal L) L(\gamma)} \mu_x(Cyl_y(\gamma^{-1}))\big
  (\sum_{f\in b(\gamma)}
  \mu_y(Cyl_y(f))\big )=\\
  \text{ by formulae in Remark~\ref{simp}}\\
  e^{-h(\mathcal L) L(\gamma)}\mu_x(Cyl_y(\gamma^{-1}))\big
  (\sum_{f\in b(\gamma)}
  e^{-h(\mathcal L) L(f)}\mu_{t(f)}(Cyl_y(f))\big )=\\
  \sum_{f\in b(\gamma)}
  e^{-h(\mathcal L) L(\gamma f)}\mu_x(Cyl_y(\gamma^{-1}))\mu_{t(f)}(Cyl_y(f))=\\
  \sum_{f\in b(\gamma)}\nu'(Cyl_X(\gamma f)).
\end{gather*}
Thus $\nu'$ is indeed a geodesic current.  We will now show that
$\nu'=\nu$, where the measure $\nu$ on $\partial^2 X$ is defined by
\thetag{$\clubsuit$}.  It suffices to show that
$\nu(Cyl_X(\gamma))=\nu'(Cyl_X(\gamma))$ for every $\gamma\in \mathcal
P(X)$.  Let $\gamma=[x,y]\in \mathcal P(X)$.

As $\nu$ is independent of the choice of $z$, we can suppose without
loss of generality that $d(z,[x,y])>0$.

Let $z'\in [x,y]$ be such that $d(z,z')=d(z,[x,y])$. Then
\begin{gather*}
  \nu(Cyl_X(\gamma))=e^{2h d(z,z')} \mu_z(Cyl_z([z,x]))
  \mu_z(Cyl_z([z,y]))=\\
  e^{2h d(z,z')} e^{-h d(z,z')}\mu_{z'}(Cyl_{z'}([z',x]))e^{-h
    d(z,z')}\mu_{z'}(Cyl_{z'}([z',y]))=\\
  \mu_{z'}(Cyl_{z'}([z',x]))\mu_{z'}(Cyl_{z'}([z',y]))=\\
  e^{-h d(z',x)}\mu_{x}(Cyl_{z'}([z',x]))e^{-h
    d(z',y)}\mu_{y}(Cyl_{z'}([z',y]))=\\
  e^{-h d(x,y)}\mu_{x}(Cyl_{z'}([z',x]))\mu_{y}(Cyl_{z'}([z',y]))=\\
  e^{-h d(x,y)}\mu_{x}(Cyl_y([y,x]))\mu_{y}(Cyl_x([x,y]))=\\
  e^{-h
    L(\gamma)}\mu_{x}(Cyl_y(\gamma^{-1}))\mu_{y}(Cyl_x(\gamma))=\nu'(Cyl_X(\gamma)).
\end{gather*}
Therefore $\nu=\nu'$, which completes the proof of
Proposition~\ref{prop:pscurr}
\end{proof}

\section{The Culler-Vogtmann outer space}

The Culler-Vogtmann outer space, introduced by Culler and Vogtmann in
a seminal paper~\cite{CV}, is a free group analogue of the
Teichm\"uller space of a closed surface of negative Euler
characteristic. We refer the reader to the original paper \cite{CV}
and to a survey paper \cite{Vog} for a detailed discussion of the
basic facts listed in this section and for the further references.

\begin{defn}[Outer space]
  Let $F$ be a free group of finite rank $k\ge 2$.  A \emph{marked
    metric graph structure} on $F$ is a pair $(\alpha, \mathcal L$),
  where $\alpha: F\to \pi_1(\Gamma,p)$ is a simplicial chart for $F$
  and $\mathcal L$ is a metric structure on $\Gamma$.  A marked
  metric graph structure is \emph{minimal} if $\Gamma$ has no
  degree-one and degree-two vertices.
  
  Two marked metric graph structures $(\alpha_1: F\to
  \pi_1(\Gamma_1,p_1), \mathcal L_1)$ and $(\alpha_2: F\to
  \pi_1(\Gamma_2,p_2), \mathcal L_2)$ are \emph{equivalent} if there
  exist an isometry $\iota: (\Gamma_1,\mathcal L_1)\to
  (\Gamma_2,\mathcal L_2)$ and a path $v$ from $\iota(p_1)$ to $p_2$
  in $\Gamma_2$ such that
\[
(\iota_{\#}\circ \alpha_1)(f)= v \alpha_2(f) v^{-1}
\]
for every $f\in F$. Clearly, minimality is preserved by equivalence of
marked metric graph structures.

The \emph{Culler-Vogtmann outer space $CV(F)$} consists of equivalence
classes of all volume-one minimal marked metric graph structures on
$F$.
\end{defn}

\begin{defn}[Elementary charts]\label{defn:ec}
  Let $\alpha: F\to \pi_1(\Gamma,p)$ be a simplicial chart for $F$,
  where $\Gamma$ has no degree-one and degree-two vertices.
  
  For each non-singular semi-metric structure $\mathcal L$ on $\Gamma$
  let $\Gamma'$, $\mathcal L'$ and $q$ be as in
  Convention~\ref{conv:semi}.  Then $q_{\#}\circ \alpha: F\to
  \pi_1(\Gamma',q(p))$ is a simplicial chart for $F$ and $(q_{\#}\circ
  \alpha, \mathcal L')$ is a minimal marked metric graph structure on
  $F$.

  Denote by $S(\Gamma)$ the set of all volume-one non-singular
  semi-metric structures on $\Gamma$.  Note that if $\Gamma$ has $N$
  non-oriented edges, then $S(\Gamma)$ is embedded as a subset of
  $\mathbb R^n$. We topologize $S(\Gamma)$ accordingly.
  
  It is not hard to see that for two non-singular semi-metric
  structures $\mathcal L_1, \mathcal L_2$ on $\Gamma$ the pairs
  $(q_{\#}\circ \alpha, \mathcal L_1')$ and $(q_{\#}\circ \alpha,
  \mathcal L_2')$ are equivalent if and only if $\mathcal L_1=\mathcal
  L_2$.  Thus $\alpha$ defines an injective map $\lambda_\alpha:
  S(\Gamma)\to CV(F)$, $\lambda_\alpha: \mathcal L\mapsto (q_{\#}\circ
  \alpha, \mathcal L')$. This map $\lambda_\alpha$ is called the
  \emph{elementary chart in $CV(F)$} corresponding to $\alpha$.
\end{defn}

Let now $S_+(\Gamma)$ denote the set of all metric structures on
$\Gamma$.  If $\Gamma$ has $n$ oriented edges then $S_+(\Gamma)$ is an
open simplex of dimension $n/2-1$ in $\mathbb R^n$ and $S_+(\Gamma)$
is dense in $S(\Gamma)$.

\begin{defn}[Topology on the outer space]
  The outer space $CV(F)$ is endowed with the weakest topology for
  which every elementary chart is a topological embedding.
\end{defn}

As explained in~\cite{CV}, the outer space $CV(F)$ is a union of open
simplices of the form $\lambda_\alpha(S_+(\Gamma))$, where
$\lambda_\alpha$ is as in Definition~\ref{defn:ec}. 

One can also view
$CV(F)$ as the space of projectivized hyperbolic length functions on
$F$ corresponding to free and discrete isometric actions of $F$ on
$\mathbb R$-trees.

\begin{defn}[Projectivized length functions]
  Let $FLen(F)$ denote the space of all hyperbolic length functions
  $\ell:F \to \mathbb R$ on $F$ corresponding to free and discrete
  isometric actions of $F$ on $\mathbb R$-trees. The space $FLen(F)$
  is endowed with the weak topology of pointwise convergence.
  
  We will say that two length functions in $FLen(F)$ are equivalent if
  they are scalar multiples of each other, and will denote by $\mathbb
  PFLen(F)$ the space of equivalence classes of elements of $FLen(F)$,
  endowed with the quotient topology. The equivalence class of
  $\ell\in FLen(F)$ is denoted $[\ell]$. For each $\ell\in FLen(F)$
  there exists a free discrete minimal isometric action of $F$ on an
  $\mathbb R$-tree $X_\ell$ such that $\ell$ is the hyperbolic length
  function for this action. Moreover, the tree $X_\ell$ and the
  corresponding action of $F$ are unique up to an equivariant
  isometry. Let $\Gamma_\ell$ denote the metric graph $X_\ell/F$.
  
  Let $FLen_1(F)$ denote the set of all $\ell\in FLen(F)$ such that
  $\Gamma_\ell$ has volume one. Note that every equivalence class
  $[\ell]\in \mathbb PFLen(F)$ has a unique representative in
  $FLen_1(F)$. For each $\ell\in FLen(F)$ the action of $F$ on
  $X_\ell$ defines an isomorphism $\alpha_\ell: F\to
  \pi_1(\Gamma_\ell,p)$, where $p\in V\Gamma_\ell$. Let $\mathcal
  L_\ell$ denote the metric structure on $\Gamma_\ell$ inherited
  from $X_\ell$. Note that the equivalence class of the marked metric
  graph structure $(\alpha_\ell, \mathcal L_\ell)$ on $F$ does not
  depend on the choice of $p$.
\end{defn}

The following statement is well-known and can be derived from results of \cite{CV}.
It shows that the outer space $CV(F)$ is homeomorphic to the spaces $FLen_1(F)$
and $\mathbb PFLen(F)$.

\begin{prop}\label{prop:cv}
\begin{enumerate}
\item The restriction of the quotient map $[\, ]: FLen(F)\to \mathbb
  PFLen(F)$ to $FLen_1(F)$ is a homeomorphism whose image is $\mathbb
  PFLen(F)$. Thus $FLen_1(F)$ is canonically homeomorphic to $\mathbb
  PFLen(F)$.
\item Let $\varrho: FLen_1(F)\to CV(F)$ be the map that takes each
  $\ell\in FLen_1(F)$ to the equivalence class of the marked structure
  $(\alpha_\ell, \mathcal L_\ell)$ on $F$. Then $\varrho: FLen_1(F)\to
  CV(F)$ is a homeomorphism.
\end{enumerate}
\end{prop}

\section{Proof of the main result}

If $(\alpha_1:F\to \pi_1(\Gamma_1,p_1), \mathcal L_1)$ and
$(\alpha_2:F\to \pi_1(\Gamma_2,p_2), \mathcal L_2)$ are two equivalent
pairs representing the same point $\eta\in CV(F)$, then $\mathbb
R$-trees $X_1=\widetilde \Gamma_1$ and $X_2=\widetilde \Gamma_2$ are
$F$-equivariantly isometric and the corresponding hyperbolic length
functions are equal.  Hence it follows from Proposition~2 of
Furman~\cite{Fur} (and it is also easy to see this directly) that the
projective $F$-Patterson-Sullivan currents corresponding to
$(\alpha_1,\mathcal L_1)$ and $(\alpha_2,\mathcal L_2)$ coincide (see
Definition \ref{def:Fcurr}). Hence the following map is well-defined:

\begin{defn}[Patterson-Sullivan map and Hausdorff dimension map]\label{def:PSHmap}
  Let $F$ be a free group of finite rank $k\ge 2$ and let $CV(F)$
  denote the outer space.

  Let $\eta\in CV(F)$. Thus $\eta$ is represented as an equivalence
  class of $(\alpha, \mathcal L)$, where $\alpha: F\to
  \pi_1(\Gamma,p)$ is a simplicial chart on $F$ such that $\Gamma$ is
  a finite connected graph without degree-one and degree-two vertices
  and where $\mathcal L$ is a volume-one metric structure on
  $\Gamma$. Consider $X=\widetilde \Gamma$ and let $d$ be the metric on $X$
  induced by $\mathcal L$.  Define $\tau(\eta)$ to be the projective
  $F$-Patterson-Sullivan current on $F$ corresponding to $(\alpha,
  \mathcal L)$. Also define $h(\eta)$ to be the Hausdorff dimension of
  $\partial X$ (which, as we have seen, is equal to the critical
  exponent $h(\mathcal L)$.)
  
  This defines a map $\tau: CV(F)\to \mathbb PCurr(F)$, which we will
  call the \emph{Patterson-Sullivan map}, and a map $h:CV(F)\to
  \mathbb R$, which we will call the \emph{Hausdorff dimension map}.
\end{defn}

\begin{thm}\label{thm:main}
  The Patterson-Sullivan map $\tau: CV(F)\to \mathbb P Curr(F)$ is a
  topological embedding.  The Hausdorff dimension map
  $h:CV(F)\to\mathbb R$ is continuous and, moreover, the restriction
  of $h$ to any open simplex in $CV(F)$ is real-analytic.
\end{thm}

\begin{proof}
  Since $CV(F)$ is locally compact, in order to prove that $\tau$ is a
  topological embedding it suffices to show that $\tau$ is continuous
  and injective.
  
  Recall the
  identification of $CV(F)$ with $\mathbb PFLen(F)$ from
  Proposition~\ref{prop:cv}.  If $\tau([\ell_1])=\tau([\ell_2])$ for
  $\ell_1,\ell_2\in FLen(F)$ then Theorem~2 of Furman \cite{Fur} implies that
  there is $r>0$ such that $r \ell_1=\ell_2$ and hence
  $[\ell_1]=[\ell_2]$, and $\tau$ is injective.
  
  We now establish that $\tau$ and $h$ are continuous.  Since every
  point of the outer space is contained in finitely many elementary
  charts, it suffices to prove that $\tau$ and $h$ are continuous on
  the image of every elementary chart in $CV(F)$.
  
  Let $\alpha: F\to \pi_1(\Gamma,p)$ be a simplicial chart for $F$,
  where $\Gamma$ has no degree-one and degree-two vertices. Let
  $\lambda_\alpha$ be the elementary chart in $CV(F)$ determined by
  $\alpha$. Recall that the image $Im(\lambda_\alpha)$ of
  $\lambda_\alpha$ consists of all points of $CV(F)$ corresponding to
  volume-one semi-metric structures on $\Gamma$ where all the edges
  with zero length are contained in a (possibly empty) subtree of
  $\Gamma$.  Corollary~\ref{cor:use} and formula $\thetag{\dag}$ in
  Proposition~\ref{prop:pscurr} imply that
  $\tau|_{Im(\lambda_\alpha)}$ and $h|_{Im(\lambda_\alpha)}$ are
  continuous and, moreover, the restriction of $h$ to the interior of
  $Im(\lambda_\alpha)$ is real-analytic.
\end{proof}

\begin{rem}[$Out(F)$-Equivariance]
  It is easy to see that the Patterson-Sullivan map $\tau$ is
  equivariant with respect to the left action of $Out(F)$ and, in
  fact, a similar statement holds in the general word-hyperbolic
  context considered by Furman~\cite{Fur}. It is even easier to see
  that $h$ is constant on each $Out(F)$-orbit and thus factors to a
  continuous map on the moduli space $\overline h:
  CV(F)/Out(F)\to\mathbb R$.

  Indeed, suppose $(\alpha: F\to \pi_1(\Gamma_,p), \mathcal L)$
  represents a point $\eta\in CV(F)$ and let $\phi\in Aut(F)$.  Let
  $X=\widetilde \Gamma$, equipped with the metric $d$ induced by
  $\mathcal L$.
  
  By definition of the left action of $Aut(F)$ (and of $Out(F)$) on
  $CV(F)$, the point $\phi \eta \in CV(F)$ is the equivalence class of
  $(\phi^{-1}\circ \alpha, \mathcal L)$. For both $\eta$ and
  $\phi\eta$ the metric graph $(\Gamma, \mathcal L)$ is the same.
  This already implies that $h(\eta)=h(\phi\eta)$.

    The action of $F$ on $X$ corresponding to $\phi\eta$ is obtained
  from the $F$-action on $X$ corresponding to $\eta$ by a
  pre-composition with $\phi^{-1}$.  The definitions imply that if
  $(\mu_x)_x$ is a family of Patterson-Sullivan measures on $\partial
  X$ corresponding to the action of $F$ on $X$ via $\alpha$, then
  $(\mu_x)_x$ is also a family of Patterson-Sullivan measures on
  $\partial X$ corresponding to the action of $F$ on $X$ via
  $\phi^{-1}\circ\alpha$.  Hence if $\nu$ is an $X$-Patterson-Sullivan
  current corresponding to the action of $F$ on $X$ via $\alpha$, then
  $\nu$ is also an $X$-Patterson-Sullivan current corresponding to the
  action of $F$ on $X$ via $\phi^{-1}\circ\alpha$.
  
  Denote $\nu_1:=\hat \alpha_{\ast}(\nu)$ and $\nu_2:=[\hat
  \phi^{-1}\circ \hat \alpha]_{\ast}(\nu)$, so that
  $\tau(\eta)=[\nu_1]$ and $\tau(\phi \eta)=[\nu_2]$. Definitions then
  imply that $\nu_2=(\hat \phi^{-1})_{\ast} \nu_1$, that is, for any
  Borel subset $A\subseteq \partial^2 F$ we have $\nu_2(S)=\nu_1(\hat
  \phi^{-1}(A))$.  By definition of the left action of $Aut(F)$ on
  $Curr(F)$ (see \cite{Ka1}) we have $(\phi\nu_1)(A)=\nu_1(\hat
  \phi^{-1}(A))$. Thus $\nu_2=\phi \nu_1$ and hence
  $\tau(\phi\eta)=\phi(\tau\eta)$, as claimed.
\end{rem}

\section{The minimal volume entropy problem}

Our goal in this section is to prove parts (1) and (3) of
Theorem~\ref{C} from the Introduction.
For the remainder of this section let $k\ge 2$, and let $\Gamma$ be
  a finite connected graph whose fundamental group $F=\pi_1(\Gamma,p)$
  with respect to a base vertex $p\in V\Gamma$ is free of rank $k$.
  Let $X=\widetilde{(\Gamma,p)}$, and let $x_0\in VX$ be a fixed lift
  of~$p$.

Let $w\in \mathcal P(\Gamma)$ and $e\in E\Gamma$. We denote by
$\langle e,w\rangle$ the number of occurrences of $e$ in $w$.  
  Let $\mathcal L$ be a positive quasi-metric structure on $\Gamma$.
  Let $w\in \mathcal P(\Gamma)$ be a reduced path.  Then
\[
L_{\mathcal L}(w)=\sum_{e\in E\Gamma} \langle e,w\rangle L_{\mathcal
  L}(e).
\]

The key step in the proof of Theorem~\ref{C} is the following
statement, which provides a sharp bound for the volume entropy of
(regular) $m$-valent metric graphs. Note that an $m$-valent graph with
the fundamental group of rank $k$ has $m(k-1)/(m-2)$ non-oriented
edges.

\begin{prop}\label{prop:ineq} For $m\ge 3$ suppose $\Gamma$ is a finite regular $m$-valent
  graph (i.e., every vertex has degree $m$) with
  fundamental group free of rank $k\ge 2$.  Let $\mathcal L$ be a
  volume-one positive quasi-metric structure on $\Gamma$.  Then
\[
h({\mathcal L})\ge \frac{m(k-1)}{m-2} \log(m-1).
\]
\end{prop}

The following lemma which we will use in the proof of Proposition \ref{prop:ineq} follows directly from the definition of the matrix
$A_{\mathcal L}(s)$.

\begin{lem}\label{lem:power}
  Let $\mathcal L$ be a positive quasi-metric structure on $\Gamma$.
  Then for any integer $t\ge 1$ and for any position $ij$ we have that
\[
[A_{\mathcal L}(s)^t]_{ij}=e^{sL(e_j)}\sum_{v} e^{-sL(v)}
\]
where the summation is taken over all reduced paths $v$ of edge-length
$t$ with the first edges $e_i$ and the last edge $e_j$.
\end{lem}

\begin{proof}[Proof of Proposition \ref{prop:ineq}]
  We consider the simple non-backtracking random walk on $\Gamma$.
  This walk can be thought of as a finite state Markov process with
  the state set $E\Gamma$ and with transition probabilities, for $e,e'\in E\Gamma$, defined as
\[
p(e,e')=\begin{cases} \frac{1}{m-1}, \quad \text{if } ee'\in \mathcal
  P(\Gamma),\\
  0, \quad \text{otherwise}\ .
\end{cases}
\]

This Markov process is irreducible since for any $e,e'\in E\Gamma$
there exists a reduced path in $\Gamma$ with initial edge $e$ and
terminal edge $e'$. The graph $\Gamma$ has $(mk-m)/(m-2)$ nonoriented
edges and $(2mk-2m)/(m-2)$ oriented edges. The uniform distribution
$\mu_0$ on $E\Gamma$, given by $\mu_0(e)=\frac{m-2}{2mk-2m}$ for every
$e\in E\Gamma$, is obviously invariant with respect to our Markov
process. Since the process is irreducible, $\mu_0$ is the only
invariant distribution on $E\Gamma$.

Let $\mu$ be the distribution on $E\Gamma$ which is uniformly
distributed on the $m$ oriented edges starting with the base-vertex
$p$. In other words, $\mu(e)=1/m$ if $o(e)=p$ and $\mu(e)=0$ if
$o(e)\ne p$.

Let $w_t=e_1,\dots,e_t$ be a trajectory of our process of length $t$.
Let $\epsilon>0$ and $e\in E\Gamma$. By Law of Large Numbers we have
\[
\lim_{t\to\infty}\mathbb P_{\mu}(|\frac{\langle
  e,w_t\rangle}{t}-\frac{m-2}{2mk-2m}|>\epsilon)=0.
\]

Let $\epsilon>0$ be arbitrary. Then there is $t_0\ge 1$ such that for
every $t\ge t_0$ and for every $e\in E\Gamma$ we have
\[
\mathbb P_{\mu}(|\frac{\langle e,w_t\rangle}{t}-\frac{m-2}{2mk-2m}|\le
\epsilon)\ge \frac{1}{2}.
\]

Denote
\begin{gather*}
  R(t,\epsilon)=\{w\in \mathcal P(\Gamma) : w \text { consists of } t
  \text { edges, and for
    every } e\in E\Gamma\\
  |\frac{\langle e,w\rangle}{t}-\frac{m-2}{2mk-2m}|\le\epsilon\}.
\end{gather*}

Thus for $t\ge t_0$
\[
\# R(t,\epsilon)\ge (1/2)\cdot m\cdot (m-1)^{t-1}\ge
(m-1)^{t-1}.\tag{$\heartsuit$}
\]

The volume of $\Gamma$ is equal to one and hence $\sum_{e\in
  E\Gamma}L_{\mathcal L}(e)=2$.

Then for every $w\in R(t,\epsilon)$ we have:
\begin{gather*}\tag{$\spadesuit$}
  L_{\mathcal L}(w)=\sum_{e\in E\Gamma} \langle e,w\rangle L_{\mathcal
    L}(e)\le\\ \sum_{e\in E\Gamma} (\frac{t(m-2)}{2mk-2m}+t\epsilon)
  L_{\mathcal L}(e)=(\frac{t(m-2)}{2mk-2m}+t\epsilon)\sum_{e\in
    E\Gamma}L_{\mathcal L}(e)=\\
  (\frac{t(m-2)}{2mk-2m}+t\epsilon)\cdot
  2=\frac{t(m-2)}{mk-m}+2t\epsilon.
\end{gather*}

Let $\displaystyle c=\min_{e\in E\Gamma} e^{sL(e)}$. Then
Lemma~\ref{lem:power} together with $(\heartsuit)$ and $(\spadesuit)$
imply that for every integer $t\ge 1$
\[
\sum_{ij} (A_{\mathcal L}(s)^t)_{ij}\ge c (m-1)^{t-1}e^{-s
  \frac{t(m-2)}{mk-m}+2t\epsilon}=\frac{c}{m-1}
[(m-1)e^{-s\frac{(m-2)}{mk-m}+2\epsilon}]^t
\]
For $s=h(\mathcal L)$ the matrix $A_\mathcal L(s)$ has spectral radius
$1$ and therefore
\[
(m-1)e^{-s\frac{(m-2)}{mk-m}+2\epsilon}\le 1
\]
and hence
\[
s\ge \frac{mk-2}{m-2}(\log(m-1) +2\epsilon).
\]
Since this is true for all $\epsilon>0$, this implies that
\[
s\ge\frac{mk-m}{m-2} \log(m-1).
\]
\end{proof}

The following easy computation shows that the bound in
Proposition~\ref{prop:ineq} is realized by the \lq\lq uniform\rq\rq\ 
volume-one metric structure, where all edges have equal lengths.

\begin{lem}\label{lem:unif}
  Let $\Gamma$ be as in Proposition~\ref{prop:ineq}.
  Let $\mathcal L_0$ be the \lq\lq uniform\rq\rq\ volume-one metric
  structure on $\Gamma$, that is $L_{\mathcal
    L_0}(e)=\frac{m-2}{mk-m}$ for every $e\in E\Gamma$.  Then
\[
h(\mathcal L_0)=\frac{m(k-1)}{m-2} \log(m-1).
\]
\end{lem}
\begin{proof}
  A direct check shows that for $s_0=\frac{m(k-1)}{m-2} \log(m-1)$ we
  have
\[A_{\mathcal L_0}(s_0)\begin{bmatrix} 1\\ 1\\ \vdots\\
  1\end{bmatrix}=\begin{bmatrix} 1\\ 1\\ \vdots\\
  1\end{bmatrix}.\] Therefore the Perron-Frobenius eigenvalue of
$A_{\mathcal L_0}(s_0)$ is equal to $1$ and hence $h(\mathcal
L_0)=s_0$, as claimed.
\end{proof}

Note that for the case $m=3$, corresponding to regular trivalent
graphs, we have $\frac{mk-m}{m-2} \log(m-1)=(3k-3)\log 2$. We are now
ready to prove part (1) of Theorem~\ref{C} from the Introduction.
\begin{thm}\label{thm:min}
  Let $F_k$ be a free group of finite rank $k\ge 2$, and let $\mathcal
  M_k:=CV(F_k)/Out(F_k)$ be the moduli space.  For the function
  $\overline h: \mathcal M_k \to\mathbb R$ we have
\[
\min \overline h= (3k-3)\log 2.
\]
This minimum is realized by any regular trivalent connected graph
$\Gamma$ with $\pi_1(\Gamma)\cong F_k$ (so that $\Gamma$ has $3k-3$
non-oriented edges), where each edge of $\Gamma$ is given length
$1/(3k-3)$.
\end{thm}

\begin{proof}
  The moduli space $\mathcal M_k$ is a union of finitely many open
  simplices, corresponding to taking all volume-one metric structures
  on all the possible minimal graphs with fundamental group free of
  rank $k$.

  Let $(\Gamma,\mathcal L)\in \mathcal M_k$. Then $(\Gamma,\mathcal
  L)$ can be approximated in $\mathcal M_k$ by trivalent metric
  graphs.  By Proposition~\ref{prop:ineq} for all of these trivalent
  graphs the volume entropy is $\ge (3k-3)\log 2$. Since $\overline h$
  is continuous on $\mathcal M_k$, it follows that $\overline
  h(\mathcal L)\ge (3k-3)\log 2$ as well. Together with
  Lemma~\ref{lem:unif} this implies the conclusion of
  Theorem~\ref{thm:min}.
\end{proof}

The following is part (3) of Theorem~\ref{C} from the Introduction.

\begin{thm}\label{thm:sup}
  Let $F_k$ be a free group of finite rank $k\ge 2$, and let $\mathcal
  M_k:=CV(F_k)/Out(F_k)$ be the moduli space. Then
\[
\sup_{\mathcal M_k} \overline{h}=\infty.
\]
\end{thm}

\begin{proof}
  First note that it suffices to prove the statement of the theorem
  for $k=2$. Indeed, suppose we know that $\sup_{\mathcal M_2}
  \overline{h}=\infty$ and let $k>2$ be arbitrary.
  
  Let $(\Gamma,\mathcal L)$ be a finite volume-one connected metric
  graph with $\pi_1(\Gamma)\cong F_2$. Let $X=\widetilde
  {(\Gamma,p)}$, where $p$ is a vertex of $\Gamma$, and let $X$ be
  endowed with a metric $d_{\mathcal L}$ induced by $\mathcal L$.
  Denote by $d$ the corresponding metric on $\partial X$.
  
  Put $\Gamma_1$ to be the wedge at $p$ of the graph $\Gamma$ with
  $k-2$ loop-edges. Let $\mathcal L_1$ be the metric structure on
  $\Gamma_1$ where each of the new loop-edges is given length
  $\frac{1}{2(k-2)}$ and where $\mathcal L_1$ restricted to $\Gamma$
  is $\mathcal L/2$. Then $\mathcal L_1$ has volume one and
  $\pi_1(\Gamma_1,p)\cong F_k$.

  Let $X_1=\widetilde{(\Gamma_1, p)}$, endowed with the induced metric
  $d_{\mathcal L_1}$. Denote by $d_1$ the corresponding metric on
  $\partial X_1$. Then $X_1$ contains an isometrically embedded copy
  of $(X,\frac{d_{\mathcal L}}{2})$ and hence $(\partial X_1,d_1)$
  contains an isometrically embedded copy of $(\partial X, d^{1/2})$.
  Taking the square root of a metric doubles the Hausdorff dimension
  and therefore
\[
h(X_1)\ge 2h(X).
\]
In particular, $h(X)\to\infty$ implies $h(X_1)\to\infty$.

Thus we may assume that $k=2$. Let $\Gamma$ be the wedge of two
loop-edges at a single vertex. Denote $E\Gamma=\{g,\bar g, f, \bar
f\}$. Let $\mathcal L$ be a volume-one metric structure on $\Gamma$
and denote $x=L(g),y=L(f)$, so that $x+y=1$ and $x,y>0$.  Then
$\overline h(\mathcal L)$ is the unique number $s>0$ such that
$\Phi_\mathcal L(s)=1$. The condition $\Phi_\mathcal L(s)=1$ is
equivalent to the existence of a positive vector $Y>0$ such that
$A_{\mathcal L}(s)Y=Y$.

The symmetry considerations imply that $Y_g=Y_{\bar g}$ and
$Y_{f}=Y_{\bar f}$. Denote $a=Y_{g}=Y_{\bar g}$ and $b=Y_{f}=Y_{\bar
  f}$. Then the system $A_{\mathcal L}(s)Y=Y$ is:

\[
\begin{cases}
  e^{-sx}(a+2b)=a,\\
  e^{-sy}(b+2a)=b,
\end{cases}
\]
Up to re-scaling we may assume $b=1$, so that the above system
transforms into the equation
\[
4=(e^{sx}-1)(e^{sy}-1).\tag{\#}
\]

Since the volume is equal to one we have $y=1-x$. For $0<x<1$ denote
by $s(x)$ the unique value $s>0$ such that the equation \thetag{\#}
holds.

We claim that $s(x)\to\infty$ as $x\to 0+$. Indeed, suppose not.  Then
there exists a sequence $x_n>0$, with $\displaystyle\lim_{n\to \infty}
x_n=0$ such that for the corresponding values $s_n=s(x_n)$ we have
$s_n\le M$, where $0<M<\infty$. Also, denote $y_n=1-x_n$. Then
$e^{s_ny_n}-1\le e^M-1=:K$. Since $0<s_n\le M$ and $\lim_{n\to \infty}
x_n=0$, we have $\lim_{n\to\infty}e^{s_nx_n}-1=0$. Therefore there
exists $m>1$ such that $0< e^{s_mx_m}-1<1/K$. Together with
$0<e^{s_my_m}-1\le K$ this implies
\[
(e^{s_mx_m}-1)(e^{s_my_m}-1)\le K\cdot (1/K)=1<4,
\]
yielding a contradiction.
\end{proof}

\section{Uniqueness of critical points for  volume entropy}

In this Section we will compute the derivative of the volume entropy
function $h$ (see also \cite{Riv} for another proof) and, as
a consequence, prove Part (2) of Theorem~\ref{C}.
For the remainder of this section, unless specified otherwise let
  $\Gamma$ be a finite connected graph without degree-one vertices and
  with the fundamental group free of rank $k\ge 2$. Let $n$ be the
  number of oriented edges of $\Gamma$ and let $N=n/2$ be the number
  of non-oriented edges of $\Gamma$. Let $Q$ be the set of all
  volume-one positive quasi-metric structures on $\Gamma$. We identify
  $Q$ with the open simplex of dimension $n-1$ in $\mathbb R^n$:
\[
Q=\{(L_1,\dots, L_n)\in \mathbb R^n: \sum_{i=1}^n L_i=2, \text{ and }
L_i>0 \text{ for } i=1,\dots, n\}.
\]

We will be using the notations from Section~\ref{sect:pf}. The proof
of Proposition~\ref{det} shows that $h$ extends to a smooth function
$s: V\to\mathbb R$ where $V$ is an open neighborhood of $Q$ in
$\mathbb R^n$. Recall the notation $s=s(\overline L)$ from
Corollary~\ref{cor:use}.

We are now going to compute the partial derivatives $\frac{\partial
  s}{\partial L_i}$. For that we will use the proof of
Proposition~\ref{det} and the following statement from Seneta's book
\cite{Se}:

\begin{prop}[see Theorem 1.5 in \cite{Se}]\label{prop:se}
  Let $A\ge 0$ be a nonzero irreducible $n\times n$ matrix. Then
\[
adj(r(A)I_n-r(A))=\epsilon YZ,
\]
where $Z$ and $Y$ are left and right Perron-Frobenius eigenvectors of
$A$ accordingly and where $\epsilon\in \{-1,1\}$.
\end{prop}

\begin{prop}\label{prop:deriv}
  Let $\mathcal L=(L_1,\dots, L_n)\in Q$. Then, at the point $\mathcal
  L$, for $i=1,\dots, n$ we have \[ \frac{\partial s}{\partial L_i}=
  sb_ie^{-sL_i}=2\epsilon sz_i||Y||^2e^{-sL_i}\det(J)^{-1},\] where
  $Z=[z_1,\dots, z_n]$ is a right Perron-Frobenius eigenvector of
  $A=A_{\mathcal L}(s)$; where $J$ and $Y$ are as in the proof of
  Proposition~\ref{det}; where $[b_1\dots b_n\ b_{n+1}]$ is the last
  row of $J^{-1}$ and where $\epsilon\in \{-1,1\}$.
\end{prop}

\begin{proof}
  Recall that in the proof of Proposition~\ref{det} the function
  $s:V\to \mathbb R$ is defined via the Implicit Function Theorem
  applied to the system of equations (!). By the Implicit Function
  Theorem via differentiating (!) we also have at $\mathcal L\in Q$:
\begin{gather*}
  \frac{\partial (y_1,\dots, y_n,s)}{\partial (L_1,\dots,
    L_n)}=-\left[ \frac{\partial(F_1,\dots, F_n,F_{n+1})}{\partial
      (y_1,\dots, y_n,s)}\right]^{-1} \frac{\partial(F_1,\dots, F_n,
    F_{n+1})}{\partial(L_1,\dots, L_n)}=\\
  =-J^{-1}
\begin{bmatrix} -se^{-sL_1} & 0 & 0 & \dots & 0\\
  0 & -se^{-sL_2}& 0& \dots & 0 \\
  \vdots & \vdots & \vdots & \dots& \vdots\\
  0 & 0 & 0 & \dots & -se^{-sL_n}\\
  0 & 0 & 0&\dots & 0
\end{bmatrix}
\end{gather*}
where \[J=\begin{bmatrix} a_{11}-1 & a_{12} & a_{13} & \dots & a_{1n}&
  -L_1y_1\\
  a_{21} & a_{22}-1& a_{23}& \dots & a_{2n} & -L_2 y_2\\
  \vdots & \vdots & \vdots & \dots& \vdots& \vdots\\
  a_{i1} & a_{i2} & a_{i3} & \dots & a_{in}& -L_iy_i\\
  \vdots & \vdots & \vdots & \dots& \vdots& \vdots\\
  a_{n1} & a_{n2} & a_{n3} & \dots & a_{nn}-1& -L_ny_n\\
  2y_1 & 2y_2 & 2y_3&\dots & 2y_n& 0
\end{bmatrix}.\]

Thus to compute $\frac{\partial s}{\partial(L_1,\dots,L_n)}$ we need
to know the last row $[b_1\dots b_n\ b_{n+1}]$ of $J^{-1}$. Then
\[
\frac{\partial s}{\partial L_i}=sb_ie^{-sL_i}\tag{+}
\] for $i=1,\dots, n$.
Note that the value $b_{n+1}$ is actually not needed since the last
row of $\frac{\partial(F_1,\dots, F_n, F_{n+1})}{\partial(L_1,\dots,
  L_n)}$ consists entirely of zeros.

Since $J^{-1}=\det(J)^{-1} adj(J)$, we have:
\[
b_i=-\det(J)^{-1} \det\begin{bmatrix} a_{11}-1 & a_{12}
  & a_{13} & \dots & a_{1n}\\
  a_{21} & a_{22}-1& a_{23}& \dots & a_{2n}\\
  \vdots & \vdots & \vdots & \dots& \vdots\\
  2y_1 & 2y_2 & 2y_3 & \dots & 2y_n\\
  \vdots & \vdots & \vdots & \dots& \vdots\\
  a_{n1} & a_{n2} & a_{n3} & \dots & a_{nn}-1
\end{bmatrix}\tag{!!}
\]
where the $\{y_j\}$ occur in the $i$-th row. Recall that $Y=[y_1, y_2\dots,
y_n]^T$ is a right Perron-Frobenius eigenvector of $A=A_{\mathcal
  L}(s)$. By Proposition~\ref{prop:se} we have $adj(A-I_n)=-\epsilon
YZ$ where $Z=[z_1,\dots, z_n]$ is a left Perron-Frobenius eigenvector
of $A$. Denote $R=adj(A-I_n)$. Then, by taking the $i$-th row
decomposition of the determinant in (!!)  we see that $b_i$ is equal
to the scalar product of $-\det(J)^{-1}[2y_1, \dots, 2y_n]$ and the
$i$-th column of $R$. Since $R=-\epsilon YZ$, the $i$-th column of $R$
is $-\epsilon z_iY$. Hence $b_i=2\epsilon z_i||Y||^2 det(J)^{-1}$.

Then by (+) we have $\frac{\partial s}{\partial
  L_i}=sb_ie^{-sL_i}=2\epsilon sz_i||Y||^2e^{-sL_i}\det(J)^{-1}$, as
required.
\end{proof}

\begin{prop}\label{prop:equal}
  Suppose that $\Gamma$ is $m$-regular for some $m\ge 3$. Let
  $\mathcal L\in Q$ be a critical point of $s|_Q$. Then
  $L_1=L_2=\dots=L_n$.
\end{prop}
\begin{proof}
  
  By the Lagrange multipliers method at a critical point of $s|_Q$ we
  have that the gradient of $s$ is parallel to the vector $(1,1,\dots,
  1)$, that is, at such a point
\[
\frac{\partial s}{\partial L_1}=\frac{\partial s}{\partial
  L_2}=\dots=\frac{\partial s}{\partial L_n}.\tag{!!!}
\]
By the earlier computations we have
$\frac{\partial s}{\partial L_i}=sb_ie^{-sL_i}=2\epsilon
sz_i||Y||^2e^{-sL_i}\det(J)^{-1}$. Then by (!!!) at a critical point
$\mathcal L$ of $s|_Q$ we have
\[
z_1e^{-sL_1}=z_2e^{-sL_2}=\dots =z_ne^{-sL_n},
\]
where $Z=[z_1,z_2,\dots, z_n]$ is a left Perron-Frobenius eigenvector
of $A$. Up to rescaling we can choose $Z$ so that the first coordinate
of $Z$ is $e^{sL_1}$. Then by the above equation
$Z=[e^{sL_1},e^{sL_2},\dots, e^{sL_n}]$.

Now we use the fact that $ZA=Z\ diag(e^{-sL_i}) M=Z$ which translates
into $[1,1,...,1]M=Z$.  Since $\Gamma$ is $m$-regular, in every column
of $M$ there are exactly $m-1$ nonzero entries, each equal to $1$.
Then $[1,1,...,1]M=Z$ translates into $z_1=\dots =z_n=m-1$. Since we
have chosen $Z$ so that $z_i=e^{sL_i}$, this implies $L_1=L_2=\dots
=L_n$, as required.
\end{proof}

The following statement is a corollary of a result obtained by
Robert~\cite{R} and Rivin~\cite{Riv}.

\begin{prop}\label{prop:robert}
  Let $\Gamma$ be a finite connected graph without degree-one vertices
  with $\pi_1(\Gamma)$ free of rank $\ge 2$ and with $N$ non-oriented
  edges.
  
  Let $U=U(\Gamma)$ be the set of all metric graph structures on
  $\Gamma$. Thus $U$ is identified with
\[
U=\{(x_1,\dots,x_N)\in \mathbb R^N: x_i>0, i=1,\dots, N\}.
\]

Then the volume entropy $h: U\to\mathbb R$ is a convex function on
$U$.
\end{prop}

The following is part (2) of Theorem~\ref{C} from the Introduction.
\begin{thm}\label{thm:unique}
  Let $\Gamma$ be a finite connected graph without degree-one and
  degree-two vertices and with fundamental group free of rank $k\ge
  1$. Let $\mathcal L$ be a volume-one metric structure such that
  $h(\mathcal L)=(3k-3)\log 2$. Then $\Gamma$ is a regular trivalent
  graph and all edges of $\Gamma$ have equal lengths in $\mathcal L$.
\end{thm}
\begin{proof}
  {\bf Case 1.} Suppose first that $\Gamma$ is trivalent. Then
  Proposition~\ref{prop:ineq} implies that $\mathcal L$ is a point of
  minimum of $s|_Q$. Therefore $\mathcal L$ is a critical point of
  $s|_Q$. Hence by Proposition~\ref{prop:equal} all edges have equal
  lengths in $\mathcal L$.
  
  {\bf Case 2.} Suppose now that $\Gamma$ is not trivalent. Then there
  exists a trivalent graph $\Gamma'$ and a sequence of volume-one
  metric structures $\mathcal L_t'$ on $\Gamma'$ such that $(\Gamma',
  \mathcal L_t')$ converges to $(\Gamma,\mathcal L)$ in the moduli
  space $CV(F_k)/Out(F_k)$ as $t\to\infty$. Moreover, since $\Gamma$
  is not trivalent, there is an edge $e$ of $\Gamma'$ whose $\mathcal
  L_t$-length converges to $0$ as $t\to\infty$, so that, after
  possibly passing to a subsequence, the points $\mathcal L_t'$
  converge to a point in the boundary of the open simplex
  $U(\Gamma')$.
  
  By continuity of $h$ we have $\displaystyle\lim_{t\to\infty}
  h(\mathcal L_t')=(3k-3)\log 2$. Let $\mathcal L_0$ be the volume-one
  metric structure on $\Gamma'$ where all edges have equal length. Let
  $S$ be the set of all points in $U=U(\Gamma)$ at distance
  $\epsilon>0$ where $\epsilon$ is smaller than the distance from
  $\mathcal L_0$ to the boundary of $U$. For sufficiently large $t$
  the segments $[\mathcal L_0, \mathcal L_t']$ intersect $S$ in one
  point denoted $\mathcal L_t''$. By convexity of $h$ we have
  $h(\mathcal L_t'')\le \max\{h(\mathcal L_0), h(\mathcal L_t')\}$.
  Since $S$ is compact, after passing to a subsequence we have
  $\displaystyle\lim_{t\to\infty} \mathcal L_t''=\mathcal L''\in S$.
  Since $\displaystyle\lim_{t\to\infty} h(\mathcal L_t')=h(\mathcal
  L_0)=(3k-3)\log 2$, it follows that $h(\mathcal L'')\le (3k-3)\log
  2$ and therefore $h(\mathcal L'')= (3k-3)\log 2$ by
  Theorem~\ref{thm:min}. By Case~1 all edges of $\mathcal L''$ have
  equal length, which contradicts the fact that $\mathcal L''\ne
  \mathcal L_0$.
\end{proof}

\end{document}